\magnification  1200\ifx\eplain\undefined \input eplain \fi

\baselineskip13pt 
\overfullrule=0pt

\def\Im{  {\rm {Im}\,}}
\def\wtilde{\widetilde}

\def\bar{\overline}
\font\smalsmalbf=cmbx8
                                                                                                                                                                                                                                                                                                                                                                                                                                                                                                                                                                                                                                                                                                                                                                                                                                                                                                                                                                                                                                                                                                                                                                                                                                                                                                                                                                                                                                                                                                                                                                                                                                                                                                                                                                                                                                                                                                                                                                                                                                                                                                                                                                                                                                                                                                                                                                                                                                                                                                                                                                                                                                                                                                                                                                                                                                                                                                                                                                                                                                                                                                                                                                                                                                                                                                                                                                                                                                                                                                                                                                                                                                                                                                                                                                                                                                                                                                                                                                                                                                                                                                                                                                                                                                                                                                                                                                                                                                                                                                                                                                                                                                                                                                                                                                                                                                                                                                                                                                                                                                                                                                                                                                                                                                                                                                                                                                                                                                                                                                                                                                                                                                                                                                                                                                                                                                                                                                                                                                                                                                                                                                                                                                                                                                                                                                                                                                                                                                                                                                                                                                                                                                                                                                                                                                                                                                                                                                                                                                                                                                                                                                                                                                                                                                                                                                                                                                                                                                                                                                                                                                                                                                                                                                                                                                                                                                                                                                                                                                                                                                                                                                                                                                                                                                                                                                                                                                                                                                                                                                                                                                                                                                                                                                                                                                                                                                                                                                                                                                                                                                                                                                                                                                                                                                                                                                                                                                                                                                                                                                                                                                                                                                                                                                                                                                                                                                                                                                                                                                                                                                                                                                                                                                                                                                                                                           
\def\z{\zeta}


\font\smalltenrm=cmr8
\font\smallteni=cmmi8
\font\smalltensy=cmsy8
\font\smallsevrm=cmr6   \font\smallfivrm=cmr5
\font\smallsevi=cmmi6   \font\smallfivi=cmmi5
\font\smallsevsy=cmsy6  \font\smallfivsy=cmsy5
\font\smallsl=cmsl8      \font\smallit=cmti8

\def\smallfonts{\lineadj{80}\textfont0=\smalltenrm  \scriptfont0=\smallsevrm
                \scriptscriptfont0=\smallfivrm
    \textfont1=\smallteni  \scriptfont1=\smallsevi
                \scriptscriptfont0=\smallfivi
     \textfont2=\smalltensy  \scriptfont2=\smallsevsy
                \scriptscriptfont2=\smallfivsy
      \let\it\smallit\let\sl\smallsl\smalltenrm}

\def\en{\eta_{n+1}}

\font\smathbold=msbm5\font\mmathbold=msbm8
\def\P{{\bf P}}\def\g{{\bf g}}\def\X{{\bf X}}
\def\Ca{{\cal C}}\font\mathbold=msbm9 at 10pt\def\Hn{{\hbox{\mathbold\char72}}^n}\def\Sn{{S^{2n+1}}}
\font\smathbold=msbm7
\def \C{{\hbox{\mathbold\char67}}}
\def\N{{\hbox{\mathbold\char78}}}\def\S{{\cal S} }
\def\R{{\hbox{\mathbold\char82}}}
\def\sR{{\hbox{\smathbold\char82}}}\def\mR{{\hbox{\mmathbold\char82}}}

\def\zbe{\z\!\cdot\bar\eta}
\def\avg{-\hskip-1.1em\int}
\def\zn{\z_{n+1}}

\def\e{\epsilon}

\font\gothic=eusm8
\def\n{{\hbox{\gothic\char78}}}

 \def\g{{\bf g}}


\def\lineadj#1{\normalbaselines\multiply\lineskip#1\divide\lineskip100
\multiply\baselineskip#1\divide\baselineskip100
\multiply\lineskiplimit#1\divide\lineskiplimit100}

\def\remark#1.{\medskip{\noin\bf Remark #1.\enspace}}
\def\endpf{$$\eqno/\!/\!/$$}

\def\pf#1.{\smallskip\noin{\bf  #1.\enspace}}

\def\noin{\noindent}
\def\lap{\Delta}

\def\ds{\displaystyle}
\def\ts{\textstyle}

\def\e{\epsilon}\def\part{\partial_t}

\def\isn{\int_{S^n}}\def\wtilde{\widetilde}
\def\what{\widehat}
\def\Rn{\R^n}
\def\im{\int_M}
\def\irn{\int_{\sR^n}}
\def\ref#1{{\bf{[#1]}}}

\def\essup{{\rm {ess\hskip.2em sup}}}
\def\p{\partial}
\def\omegan{\omega_{2n+1}}

\def\ints{\int_\Sigma}

\def\K{{\cal K}}\def\A{{\bf A}}\def\D{{\cal D}}

 \def\F{{\cal F}}
\def\km{\wtilde k_m}

\def\op{{\hbox { Op}}}

\def\ker{{\rm Ker}}



\centerline{\bf{  Adams inequalities on measure spaces}}
\bigskip\centerline{Luigi Fontana, Carlo Morpurgo}
\vskip1em
\midinsert
{\smalsmalbf Abstract. }{\smallfonts In 1988 Adams obtained sharp Moser-Trudinger inequalities on bounded domains of $\mR^n$. The main step was a sharp exponential integral inequality for convolutions with the Riesz potential. In this paper we extend and improve  Adams' results to  functions defined on arbitrary measure spaces with finite measure. The Riesz fractional integral is replaced by  general integral operators, whose kernels satisfy suitable and explicit growth conditions, given in terms of their distribution functions; natural conditions for sharpness are also given. Most of the known results about Moser-Trudinger inequalities can be easily adapted to our unified scheme. We give some new applications of our theorems, including:   sharp higher order Moser-Trudinger trace inequalities, sharp Adams/Moser-Trudinger inequalities for general elliptic differential operators (scalar and vector-valued), for  sums of weighted potentials, and for operators  in the CR setting. 
}

\endinsert\bigskip
\centerline{\bf INTRODUCTION}\bigskip\bigskip

Exponential integrability can often compensate for  lack of boundedness, as a natural (although  weaker) condition. There are 
numerous  important instances of this idea in the literature, the first is perhaps due to Zygmund. It is well known that the conjugate function of a bounded function on the torus $T$ need not  be bounded, but in 1929 Zygmund proved that for all $\lambda < {{\pi}\over{2}}$ the conjugate function $\tilde f$  satisfies
$$
\int_{{\bf T}} \exp\Big({ \lambda |\tilde f(\theta)|}\Big) d\theta \leq C_{\lambda}
$$
whenever $f$ is real-valued and  belongs to the closed unit ball of $L^{\infty}({\bf T})$ ([Z], Ch. VII). Cancellation, through Cauchy's integral formula, plays the central role in the proof of this result.

On the other hand, size has the major role in  the chain of results that followed a 1967 paper by Trudinger, in which he showed that exponential integrability fills the gap in Sobolev's embedding
theorem (see also earlier versions in [Po] and  [Yu]):

\proclaim {Theorem [Tr]}. Let $\Omega $ be an open and  bounded set in ${\R}^n$, $n>1$.  
There exist constants $\lambda $ and $C$ such that, if $u$ belongs to the Sobolev space $W_0^{1,n}(\Omega)$ and 
$(\int_{\Omega}|\nabla u|^{n})^{{1\over n}} \leq 1$, then 
$$
\int_{\Omega} \exp\Big({\lambda |u|^{{{n}\over{n-1}}}}\Big) dx \leq C.\eqno(1)
$$
\par
In 1971 Moser sharpened the result by showing that
$
\lambda = n \omega _{n-1}^{{{1}\over{n-1}}}
$
is  best possible in (1), where  $\omega_{n-1}$
is the surface measure of the unit sphere in ${\R}^{n}$.
\eject
Due to the wide range of applications in PDE's, Differential Geometry  and String Theory, Moser's result triggered an enormous amount of work in the years that followed, and  up to present time. Several aspects and extensions of Moser's inequality were studied, and   still are part of an active field of research: existence of extremals, Neumann conditions rather than Dirichlet,     
settings other than $\Rn$, higher order derivatives, and more. All but a handful of the references listed in the back of this article 
deal with Moser-Trudinger inequalities, in some form or another, and the list is only partial.


Adams' paper in 1988, however, represents a true turning point. Not only did he extend
Moser's sharp result to higher order derivatives, but he also set the strategy that opened the way to most of the later work in the field. We recall here the basic developments. Adams' generalization of
Moser's theorem is:

\bigskip

\proclaim Theorem [Ad1].
Let $\Omega $ be  an open and bounded set   in ${\R}^{n}$, $\,n>1$, and let  $m\in \N$ with $m<n$.  There are constants $\beta(m,n)$ and $C$ with the following
property: If $u$ belongs to the Sobolev space $W_0^{m,n/m}(\Omega)$ and 
$\parallel \nabla^{m}u \parallel_{n/m}  \leq 1$, then
$$
\int_{\Omega} \exp\Big[{\beta(m,n) |u(x)|^{{{n}\over{n-m}}}}\Big] dx \leq C.\eqno(2)
$$
\par
The constant $\beta(m,n)$ is given explicitly in [Ad1] and it  is  sharp (see also Theorem 6 in Section 5). Also, $\nabla^{m}$ means $\Delta^{{m \over 2}}$
when $m$ is even and $\nabla \Delta^{{{m-1}\over {2}}}$ when $m$ is odd, where $\Delta$ denotes the positive Laplacian on $\Rn$.
              \smallskip                                                   
Adams' approach consists of five main steps. \smallskip
\noin{\bf Step 1.} Represent  $u$ in terms of $\nabla^m u$,  via  convolutions with the Riesz potential. 
\smallskip\noin{\bf Step 2.} Formulate the following sharp theorem on exponential integrability for Riesz potentials (a dual, but more general, version of the above theorem that has its own relevance).
The first theorem follows immediately from the second, apart for some
extra work necessary to ensure that the inequality is indeed sharp.

\bigskip

\proclaim  Theorem [Ad1].
For $1 < p < \infty$ , there is a constant $C$ such that for all 
$f \in L^{p} ({\R}^{n})$ with support contained in $\Omega $ and $\parallel f \parallel_{p} \leq 1 $, 
$$
\int_{\Omega} \exp\Big[{{{n}\over {\omega_{n-1}}}|I_{\alpha} \ast f(x)|^{p'} } 
 \Big]dx \leq C\eqno(3)
$$
where $\alpha = n/p $, $1/p + 1/p' = 1$, and 
$I_{\alpha} \ast f(x) = \int | x-y |^{\alpha-n} f(x) dy$.
The constant $n/\omega_{n-1}$ cannot be replaced by any larger number without  
forcing  $C$  to depend on $f$ as well as on $p$ and $n$.
\par 
\eject 
\smallskip\noin{\bf Step 3.}  The third  step of Adam's strategy  is to reduce the proof of the above theorem to a one-dimensional exponential inequality by using  a lemma due to O'Neil: if $T$ is a convolution operator on a measure space, then 
$$T(f,g)^{**}(t)\le t f^{**}(t)g^{**}(t)+\int_t^\infty f^*(u)g^*(u)du,\qquad t>0 \eqno(4)$$
where $f^*$ denotes nonincreasing rearrangement on the half-line, and $f^{**}(t)=t^{-1}\ds\int_0^t f^*(u) du$.

\bigskip
\smallskip\noin{\bf Step 4.} The next step is to prove the one-dimensional exponential inequality derived in  step 3 by means of a technical lemma,  now known as the ``Adams-Garsia's lemma". 
\medskip\noin{\bf Step 5.} The final step is to show that the exponential constant is sharp, by showing that for any larger constant one can find a suitable sequence of functions that makes the exponential integral arbitrarily large.
  \medskip
In his PhD thesis (1991) Fontana adapted Adams' strategy, and extended his results  in the setting of compact Riemannian manifolds [F]. In that situation the Green function  replaces the Riesz potential in step 1; the corresponding integral representation is no longer a global convolution,  but locally the Green kernel is a  perturbed Riesz potential. These facts eventually lead to  suitable versions of  O'Neil's lemma  and Adams-Garsia's lemma; these  modified lemmas  could not be deduced from the original ones, even though the original proofs were successfully adapted to the perturbed setting [F].

 Several other authors also  used Adams strategy, sometimes partially,  in order to prove sharp Moser-Trudinger estimates in various  settings. In most cases, like in [F], some    individual steps had to be adapted, and sometimes their proofs were only sketched, or even omitted.

Recently ([BFM]), the authors of this paper, in joint project with Tom Branson, needed a sharp form of various Moser-Trudinger inequalities in the CR setting in order to obtain the sharp version of Beckner-Onofri's inequality on the complex  sphere. Independently Cohn and Lu [CoLu 1,2] had  worked out Adams and Moser-Trudinger sharp estimates
in some very special cases, which were not suitable to our needs. While working out yet another version of Adams strategy, we realized that steps 2,3,4,5 could be formulated in an arbitrary measure space, for integral operators  more general than convolutions, and with kernels  satisfying suitable growth  and integral conditions.

 It was then that we seriously looked into the possibility of a general result that would encompass and unify the various  Adams-type  procedures,  with the hope that it would prove to be useful to authors in need of such sharp estimates in a variety of situations. Stripped down to its essence, the present paper could be summarized as follows. 

\eject

Suppose that $T$ is an integral operator of type
$$Tf(x)=\int_M k(x,y)f(y)d\mu(y)\,,\qquad x\in N$$
where $(M,\mu)$,$\,(N,\nu)$ are  measure spaces with finite measure, and suppose that the kernel $k(x,y) $ satisfies
$$\sup_{x\in N}\mu\Big(\{y\in M: \,|k(x,y)|>s\}\Big)\le A s^{-\beta}\Big(1+O(\log^{-\gamma} s)\Big)\eqno(5)$$
$$\sup_{y\in M}\nu\Big(\{x\in N:\, |k(x,y)|>s\}\Big)\le B s^{-\beta_0}\eqno(6)$$
as $s\to+\infty$, where $\beta>1$, $\beta'$ is the conjugate exponent, $0<\beta_0\le \beta$, and $\gamma>1$.
 Then we have an exponential inequality
of type 
$$\int_N \exp\bigg[{\beta_0\over A\beta} \bigg({|Tf(x)|\over\|f\|_{\beta'}}\bigg)^{\beta}\bigg]\,d\nu(x)\le C\eqno(7)$$
for any $f\in L^{\beta'}(M)$.
As for the sharpness statement, if equality holds in (5) then the constant $\beta_0/(A\beta)$ in (7) is sharp, provided that 
certain reasonable ``regularity" conditions are satisfied by the kernel $k$.

 The main feature of this result, which is Theorem 1 in the next section, is   that it reduces Moser-Trudinger inequalities 
for integral operators, or in ``dual form", to a couple of  estimates for the distribution functions of their kernels, and the sharpness result (under suitable but reasonable geometric conditions) to a single integral estimate (see d) in Theorem 4). 
In some cases estimates (5) and (6) are rather trivial to check, like for the Riesz potential $k(x,y)=|x-y|^{d-n}$, on a domain $\Omega$, for which
$$\sup_{x\in\Omega}| \{y\in\Omega :\,|x-y|^{d-n}>s\}|={\omega_{n-1}\over n}\,s^{-{n\over n-d}}.\eqno(8)$$
In other situations the  asymptotics of the distribution function of $k$ could be a bit more involved, but they are usually  a consequence of an  asymptotic expansion of the kernel $k$ around its singularity. For example, kernels that satisfy (5) and (6) are those of type
$$k(x,y)=c(d,n)|x-y|^{d-n}+O(|x-y|^{d-n+\e})\eqno(9)$$
some $\e>0$, or more generally of type  
$$k(x,y)=k_{d-n}(x,x-y)+O(|x-y|^{d-n+\e})\eqno(10)$$
some suitable $k_{d-n}(x,z)$ homogeneous of order $d-n$ in $z$ (see Lemma 9).
These are in fact more than just examples. It was already shown  by Fontana in [F]  that  one can still set up the Adams machinery for powers of Laplace-Beltrami operators on compact manifolds without boundary,  even though such operators have  fundamental solutions  that  do not satisfy the precise identity  (8), but instead satisfy a perturbed version like (9), in local coordinates. 
 
The fact that error terms are allowed in the asymptotics of the kernels or their distribution functions is an important point of our theory.   Indeed,   (10) is precisely the type of expansion satisfied by 
the classical parametrix   of elliptic pseudodifferential operators of order $d$  on bounded domains of $\Rn$ (or on compact manifolds, in local coordinates). Whenever an elliptic operator $P$, say on a domain $\Omega$, has a fundamental solution $T$ with such kernel, we can  write any compactly supported smooth function as $u=T(Pu)$ and almost  immediately obtain a sharp Moser-Trudinger inequality of type 
$$\int_\Omega \exp\bigg[A^{-1} \bigg({|u(x)|\over\|Pu\|_{p}}\bigg)^{p'}\bigg]\,dx\le C$$
where the sharp constant $A^{-1}$ depends on the principal symbol of $P$. This is in fact one of the applications we give of our main theorem, extending Adams original result (2) to  a wide class of  scalar and vector-valued elliptic differential operators (see Theorems 10, and 12). In the special case of second order elliptic operators, the sharp constant is even more explicitly described in terms of the  matrix  formed by the second order coefficients (see Corollary 11).

Another feature of our main theorem  is that it offers ample flexibility in the choice of the base measure spaces $(M,\mu)$ and $(N,\nu)$. To illustrate this point we offer an  extension of a very recent result  of Cianchi [Ci1] who proved that if $\nu$ is a Borel measure on $\Omega\subseteq \R^n$ satisfying $\nu\big(B(x,r)\cap \Omega\big)\le C r^\lambda$, for suitable $\lambda\in (0,n]$, and for small $r$, then 
$$\int_\Omega \exp\bigg[\lambda \omega_{n-1}^{1\over n-1} \bigg({|u(x)|\over\|\nabla u\|_{n}}\bigg)^{n'}\bigg]\,d\nu(x)\le C\eqno(11)$$
for all $u\in W_0^{1,n}(\Omega)$. As Cianchi observed, this result immediately leads to inequality for traces of functions, either   on boundaries of smooth $\lambda$-dimensional submanifolds of $\Rn$ or on sets of fractal type.  Cianchi's proof of the above inequality did not follow the representation formula as in Adams' original paper, step 1. By  use of a  trace Sobolev inequality also due to Adams (see (59)) and clever rearrangement results, Cianchi is however able to  make some contact with Adams' original steps 2,3,4,5.

In Theorems 6 and 7, we extend (11) to higher order operators and potentials. We especially hope to show how (11) and its higher order versions  are  part of the same large family of Adams/Moser-Trudinger inequalities, and are  in fact simple applications of our main theorems. The role of the constant is clearly explained in terms of the interactions between the base measures $d\nu(x), \,d\mu(y)=dy, $ and the Riesz potentials, as given in (5) and (6).

We would like to point out that our original formulation of (5)-(7) had $\beta=\beta_0$, and was based (among many other things)  on an improved version of O'Neil's lemma given as in (20). It was only after being aware of Cianchi's result that we started looking for a further improvement of O'Neil's lemma and (5)-(7). In particular it was Cianchi's  idea to exploit Adams' trace inequality (59) that eventually lead us to exploit instead Adams' weak-type estimates (21), in order to obtain a further substantial extension of O'Neil's lemma.

\medskip
In a third application, we consider Adams inequalities for sums of weighted Riesz potentials, i.e. for integral operators with kernel
$$K(x,y)=\sum_{j=1}^N g_j(x,y)|x+a_j-y|^{ d-n}$$
where the functions $g_j$ are H\"older continuous, and where $x$ and $y$ are allowed to move in different domains. The sharp constant is explicitly described even in this case, see Theorem~15.
\smallskip
Finally, and this was the original motivation for our work, we turn to the CR setting, by proving a sharpness result for some 
Adams' inequalities on the complex sphere, which were only partially proved  [BFM], using the methods of this work.

\medskip

The paper is organized in two main parts. In Part I we give the main results, in a measure-theoretical setting. Some portions of some  proofs are of course based on Adams' and O'Neil's original arguments,  but we decided to include them, in part because the  modifications are many, and often  not trivial,   and in part to achieve a rather self-contained and cleaner presentation. 

 In Part II we give several new applications of the general results of Part I:  higher order Adams and Moser-Trudinger trace inequalities, Moser-Trudinger and  Adams inequalities for general and then specific  elliptic operators and parametrix-like potentials respectively, followed by those for sums of weighted Riesz potentials, and finally for certain types of potentials arising in CR geometry. 
Some  of these applications  could  be combined together, but we decided to keep them separate in order to highlight the relevant aspects of a given setting, rather than presenting more comprehensive theorems with too many parameters. 

We certainly do not claim to have covered every possible Moser-Trudinger inequality, in fact we hope that many more could be  obtained using our setup, in a relatively painless way, and in a variety of settings. Moreover, in Theorems 10 and 12, a general form of the sharp exponential constant is given, but in specific cases it could be more helpful  to know this constant more explicitly. In this work we limit ourselves to give more explicit values in the case of  second order operators and certain other vector-valued inequalities, but more such computations are   possible. 
\eject
Another interesting situation arises regarding  Moser-Trudinger inequalities  in the space $W^{d,p}(\Omega)$, i.e. without boundary conditions. In [Ci2] Cianchi obtained a sharp inequality for the case $W^{1,n}(\Omega)$, but using different tools than ours, such as the isoperimetric inequality. It is possible that our methods are suitable to handle at least some special cases, such as low order operators, or particular domains.

\bigskip\centerline{\bf PART I: ABSTRACT THEOREMS}\bigskip

\noin {\bf 1.  Adams inequalities on measure spaces}\bigskip
 Let $(M,\mu)$ be a measure space, and $\mu$ a finite measure. 
Given a  measurable $f:M\to [-\infty,\infty]$ its distribution function will be denoted by 
$$m(f,s)=\mu\big(\{x\in M: |f(x)|>s\}\big),\qquad s\ge0$$
its  nonincreasing rearrangement by
$$f^*(t)=\inf\big\{s\ge0:\,m(f,s)\le t\big\},\qquad t>0$$
and 
$$f^{**}(t)={1\over t}\int_0^t f^*(s)ds,\qquad t>0$$
Given another finite measure space $(N,\nu)$ and  a $\nu\times\mu-$measurable function $k:N\times M\to[-\infty,\infty]$ we let, for $t>0$, 
$$ k_1^*(t)=\sup_{x\in N} k^*(x,\cdot)(t)$$
$$k_2^*(t)=\sup_{y\in M} k^*(\cdot,y)(t)$$
where $k^*(x,\cdot)(t)$ is the nonincreasing rearrangement of $k(x,y)$ with respect to the variable $y$ for fixed $x$, and $k^*(\cdot,y)(t)$ is its analogue for fixed $y$. With a slight abuse of notation we set
$$k_j^{**}(t)={1\over t}\int_0^t k_j^*(s)ds,\qquad t>0,\;j=1,2$$
If $k_2^*\in L^1\big([0,\infty)\big)$, or equivalently $m(k_2^*,\cdot)\in L^1\big([0,\infty)\big)$, then the integral operator
$$ Tf(x)=\int_M k(x,y)f(y)d\mu(y)\eqno(12)$$
is well defined and continuous from $L^1(M,\mu)$ to $L^1(N,\nu)$. In fact, as we shall see later, $Tf$ is also well defined on some $L^p$ under weaker integrability conditions on $k_2^*$, but with additional restrictions on $k_1^*$.

\medskip
Here is our main theorem:

\proclaim Theorem 1. Let $k:N\times M\to[-\infty,\infty]$ be measurable on the finite measure space $(N\times M,\nu\times\mu)$ 
and  such that
$$ m(k_1^*,s)\le A s^{-\beta}\Big(1+O(\log^{-\gamma} s)\Big)\eqno(13)$$
$$ m(k_2^*,s)\le B s^{-\beta_0}\eqno(14)$$
as $s\to+\infty$,  for some $\beta,\gamma>1$, $\,0<\beta_0\le \beta$ and $A,B>0$. Then, $T$ is  defined by (12) on $L^{\beta'}(M)$ and 
there exists a constant $C$ such that 
$$\int_N \exp\bigg[{\beta_0\over  A\beta}\bigg({|Tf|\over \|f\|_{\beta'}}\bigg)^\beta\,\bigg]d\nu\le C\eqno(15)$$
for each $f\in L^{\beta'}(M)$, with $\ds{1\over \beta}+{1\over\beta'}=1.$
\par

\noin{\bf Remarks.}

\noindent{\bf 1.} 
 It is possible to modify slightly the arguments in order to include in Theorem 1 the case of Lorentz spaces. For simplicity we just treat $L^p$ spaces.

\medskip
\noin{\bf 2.} Theorem 1 holds verbatim in case $k$ is complex-valued and $T$ acts on complex-valued functions, provided that $|k(x,y)|$ satisfies conditions (13), (14).\smallskip

Theorem 1, as an immediate corollary of itself, 
  can be extended to vector-valued functions as follows. For a measurable $F:M\to\R^n$, $F=(F_1,...,F_n)$, define $|F|=(F_1^2+...+F_n^2)^{1/2}$ and say $F\in L^p(M)$ if $\int_M|F|^p<\infty$, likewise for vector-valued functions defined on $N$, valued on $\R^n$, or  on $\bar \R^n=[-\infty,\infty]^n$.

\smallskip
\proclaim Theorem 1'.  Let  $K:N\times M\to\bar\R^n$, where $K=(K_1,...,K_n)$, be measurable and such that   $k(x,y)=|K(x,y)|$ satisfies conditions (13) and (14) of Theorem 1. If  
$$TF(x)=\int_M K(x,y)\cdot F(y) \,d\mu(y)=\int_M \sum_{j=1}^n K_j(x,y)F_j(y)\,d\mu(y)$$
then, $T$ is defined on $L^{\beta'}(M)$ and 
there exists a constant $C$ such that 
$$\int_N \exp\bigg[{\beta_0\over  A\beta}\bigg({|TF|\over \|F\|_{\beta'}}\bigg)^\beta\,\bigg]d\nu\le C\eqno(16)$$
for each $F\in L^{\beta'}(M)$, with $\ds{1\over \beta}+{1\over\beta'}=1.$
\par
\eject
The formulation in terms of vector-valued function is useful since in many cases one has a  representation formula
of a function which involves the  gradient operator, as in the classical Adams setting. Needless to say a similar version of the inequality holds for $\C^n$-valued kernels and functions. It is important to point out that while the inequality of Theorem~1' is an immediate consequence of the scalar case, via Cauchy-Schwarz,  this is not the case for the  sharpness statement  (see Theorem 4).
 \medskip

The following elementary facts about rearrangements will be useful ($f,g$ denote two measurable functions on $M$):
\smallskip
\noin{\bf Fact 1.}  $\;m(f,s)=m(f^*,s)$ and $\;m(f^*,s)\le m(g^*,s)$ for all $s>s_0$ (some $s_0>0$) if and only if $f^*(t)\le g^*(t)$ for all $t<t_0$ (some $t_0>0$).\smallskip
\noin{\bf Fact 2.} If $\psi(s)$ is continuous and strictly decreasing on $[s_0,\infty)$ then 
$\inf\{s:\,\psi(s)\le t\}=\psi^{-1}(t)$ for $t<\psi(s_0)$, (and hence $\psi$ is the distribution function of $\psi^{-1}$ on that interval). 
\smallskip
\noin{\bf Fact 3.} Given a measurable $k(x,y)$ on $N\times M$, if $\wtilde m(k,s)=\sup_x m\big(k^*(x,\cdot),s\big)=\sup_x m\big(k(x,\cdot),s\big)$ and $\wtilde k(t)=\inf\big\{s:\, \wtilde m(k,s)\le t\big\}$, then $m(\wtilde k,s)=\wtilde m(k,s)$ and $\wtilde k(t)=\sup_x k^*(x,\cdot)(t)$.
\smallskip
\noin{\bf Fact 4.} The following are equivalent ($A,\beta,\gamma>0$): \smallskip
\item{a)} $\;m(f^*,s)\le As^{-\beta}(1+C\log^{-\gamma}s),$ for all  $s>s_0>1$
\smallskip\item{b)} $\; f^*(t)\le A^{1/\beta} t^{-1/\beta}(1+C'|\log t|^{-\gamma}),$ for all $t<t_0<1$. 
\smallskip Likewise, the following are equivalent:\smallskip
\item{a')} $\;m(f^*,s)\ge As^{-\beta}(1-C\log^{-\gamma}s)>0,$ for all $s>s_0>1$
\smallskip\item{b')} $\; f^*(t)\ge A^{1/\beta} t^{-1/\beta}(1-C'|\log t|^{-\gamma})>0,$ for all $t<t_0<1$. 

\def\supp{{\rm supp}}

\medskip\medskip
The first, and crucial, step in the proof of Theorem 1 is the following $L^p$  version of O'Neil's lemma:

\medskip
\proclaim Lemma 2 (Improved O'Neil's Lemma). Let $k: N\times M\to [-\infty,\infty]$ be measurable and
$$k_1^*(t)\le M t^{-1/\beta},\qquad k_2^*(t)\le B t^{-1/\beta_0},\qquad t>0\eqno(17)$$
with $\beta>1$ and $0<\beta_0\le \beta$.
If  
$$\max\Big\{1,{\beta-\beta_0\over\beta-1}\Big\}< p<{\beta\over\beta-1}=\beta',\qquad\quad q={p\beta_0\over \beta-(\beta-1)p}>p\eqno(18)$$ then $T$ is defined on $L^{\beta'}(M)$, in fact $T:L^p(M)\to L^{q,\infty}(N)$ and bounded, and  there is a constant $C=C(M,B,\beta,\beta_0,p)$ such that for any $f\in L^{\beta'}(M)$
$$ (Tf)^{**}(t)\le C\,\max\big\{\tau^{-{\beta_0\over  q \beta}}, t^{-{1\over q}}\big\}\int_0^\tau f^*(u) u^{-1+{1\over p}}du+\int_\tau^\infty f^*(u)k_1^*(u)du,\quad \forall t,\tau>0.\eqno(19)$$
If instead of (17) we assume $\,k_1^*, k_2^*\in L^1\big([0,\infty)\big)$ then for every $f\in L^1(M)$
$$(Tf)^{**}(t)\le \tau\, \max\big\{k_1^{**}(\tau),k_2^{**}(t)\big\}\, f^{**}(\tau)+\int_\tau^\infty f^*(u)k_1^*(u)du,\quad \forall t,\tau>0.\eqno(20)$$
\par \smallskip
 We observe that  inequality (20) implies (19) in case $\beta_0>1$, that is when both    
 $k_1^*$ and $k_2^*$ are integrable, and it is also perfectly suitable to prove Theorem 1 in that case, but it is useless when $\beta_0\le 1$.
\medskip
\pf Proof.  We begin right away with the following weak-type estimate due to Adams [Ad3]. If $k$ and $f$ are nonnegative,  with 
$$\sup_{x\in N} \,m\big(k(x,\cdot),s)\le M s^{-\beta},\qquad \sup_{y\in M} \,m\big(k(\cdot,y),s\big)\le Bs^{-\beta_0}$$
 which are equivalent to (17), and under the hypothesis  (18), then for $s>0$
$$ s\, m(Tf,s)^{1\over q}=s\,\nu\big(\{x: Tf(x)>s\}\big)^{1\over q}\le {q^2\over \beta_0(q-p)}M^{1-{1\over p}}B^{1\over q} \|f\|_p\eqno(21)$$
or
$$ (Tf)^*(t)\le C t^{-{1\over q}} \|f\|_p,\qquad \forall t>0.\eqno(22)$$
This means that $T:L^p(M)\to L^{q,\infty}(N)$ is bounded, in particular  $T$ is well defined on $L^{\beta'}(M)\subseteq L^p(M)$. 

Without loss of generality we can assume throughout this proof that both $k$ and $f$ are nonnegative. With a slight abuse of language we let $\supp(f)=\{x\in M:\,f(x)\neq0\}$. The main step of the proof relies on the following:
\medskip
\proclaim Claim (See also Lemma 1.4 in [ON]). If $\mu(\supp f)=z$ and $0\le f(z)\le \alpha$, and if $k_1^*,k_2^*$ satisfy conditions (17), then $\forall t>0$
$$(Tf)^{**}(t)\le\alpha \,z\,k_1^{**}(z).\eqno(23)$$
$$(Tf)^{**}(t)\le C\, \alpha\, z^{{1\over p}} t^{-{1\over q}}.\eqno(24)$$
If instead of $(17)$ we assume that $k_1^*$ and $k_2^*$ are integrable, then (23) holds and (24) can be replaced by 
$$(Tf)^{**}(t)\le\alpha \,z\,k_2^{**}(t)\eqno(25).$$\par\medskip
Assuming the Claim, the proof of the lemma proceeds as follows.
For fixed $t,\tau>0$, pick  $\{y_n\}_{-\infty}^\infty$ such that $y_0=f^*(\tau),\,y_n\le y_{n+1}, \,y_n\to+\infty$
as $n\to+\infty$, and $y_n\to0$ as $n\to-\infty$.  Then 
$$f(y)=\sum_{-\infty}^\infty f_n(y)\quad{\hbox{where}}\quad f_n(y)=\cases{0 &if $\;f(y)\le y_{n-1}$\cr f(y)-y_{n-1} 
&if $\;y_{n-1}<f(y)\le y_n$ \cr y_n-y_{n-1} & if $\;y_n<f(y).$\cr}$$

Observe that supp$f_n\subseteq E_n:=\big\{y: f(y)>y_{n-1}\big\}$,  $\;\mu(E_n)=m(f,y_{n-1})$, and also  $\;0\le f_n(y)\le y_n-y_{n-1}$. Write 
$$f=\sum_{-\infty}^0f_n+\sum_1^{\infty}f_n=g_1+g_2$$
so that  $(Tf)^{**}\le (Tg_1)^{**}+(Tg_2)^{**}$ (this is the subadditivity of $(\cdot)^{**}$). Using  (24) we obtain
$$(Tg_2)^{**}(t)\le \sum_1^\infty (Tf_n)^{**}(t)\le Ct^{-{1\over q}}\,\sum_1^{\infty}(y_n-y_{n-1})\big(m(f,y_{n-1})\big)^{{1\over p}}$$
so that taking the inf over all such $\{y_n\}$ we get
$$\eqalign{&(Tg_2)^{**}(t)\le Ct^{-{1\over q}}\int_{f^*(\tau)}^\infty \big(m(f,y)\big)^{{1\over p}} dy=-\int_0^\tau\big(m(f,f^*(u))\big)^{{1\over p}}d\,f^*(u)\cr&\le- \int_0^\tau u^{{1\over p}} d\,f^*(u)=-u^{{1\over p}}f^*(u)\Big|_0^\tau+{1\over p}\int_0^\tau u^{-1+{1\over p}}f^*(u)du\le{1\over p}\int_0^\tau u^{-1+{1\over p}}f^*(u)du. \cr}  $$
(The last inequality follows since $f\in L^{\beta'}\Longrightarrow t^{{1\over \beta'}}f^*(t)\to0$, as $t\to0$.)

Likewise, using  (23)
$$(Tg_1)^{**}(t)\le \sum_{-\infty}^0 (Tf_n)^{**}(t)\le \sum_{-\infty}^0 (y_n-y_{n-1})m(f,y_{n-1})k_1^{**}\big(m(f(y_{n-1}))\big)$$
and so
$$\eqalign{&(Tg_1)^{**}(t)\le\int_0^{f^*(\tau)} m(f,y)k_1^{**}\big(m(f,y)\big)dy=
-\int_\tau^\infty m\big(f,f^*(u)\big)k_1^{**}\Big(m\big(f,f^*(u)\big)\Big) df^*(u)\cr&=-\int_\tau^\infty u \,k_1^{**}(u)df^*(u)=
-u\, k_1^{**}(u)f^*(u)\bigg|_\tau^\infty+\int_\tau^\infty k_1^*(u)f^*(u)du\cr&\le \tau\,k_1^{**}(\tau)f^*(\tau)+\int_\tau^\infty f^*(u)k_1^*(u)du\le \tau^{1-{1\over p}}\,k_1^{**}(\tau)\int_0^\tau f^*(u)u^{-1+{1\over p}}du+\int_\tau^\infty f^*(u)k_1^*(u)du\cr&
\le C\,\tau^{1-{1\over p}-{1\over \beta}}\int_0^\tau f^*(u)u^{-1+{1\over p}}du+\int_\tau^\infty f^*(u)k_1^*(u)du
}$$
and (19) follows since $\ds{{1\over p}+{1\over \beta}-1={\beta_0\over q\beta}}$.

 To prove (20), assume that $k_1^*,\,k_2^*$ and $f$ are integrable and estimate $(Tg_1)^{**}$ as before. The estimate for $(Tg_2)^{**}$ is now performed as above, but using (25) instead of (24). This yields

$$\eqalign{(Tf)^{**}(t)\le  \max\big\{k_1^{**}(\tau),k_2^{**}(t)\big\}\,\bigg[\tau f^*(\tau )&+\int_{f^*(\tau )}^\infty m(f,y)dy\bigg]+\int_\tau ^\infty f^*(u)k_1^*(u)du\cr}$$
and (20) follows from the identity
$$ \int_{f^*(\tau )}^\infty m(f,y)dy=\int_{f^*(\tau )}^\infty m(f^*,y)dy=\int_0^\tau  f^*(u)du-\tau  f^*(\tau ).$$

\medskip
\smallskip\noin{\bf Proof of Claim.}  Let $r>0$ and   set 
$$k_r(x,y)=\cases{k(x,y) & if $\;k(x,y)\le r$\cr\cr r & otherwise,\cr}\qquad k(x,y)=k_r(x,y)+k^r(x,y).$$
so that 
$$Tf(x)=\im k_r(x,y)f(y)d\mu(y)+\im k^r(x,y )f(y)d\mu(y)=h_1(x)+h_2(x).$$
Assume that $k_1^*$ is integrable. Then,
for every given $x$
$$ h_2(x)\le \|f\|_\infty^{}\im k^r(x,y)d\mu(y)\le \alpha \int_r^\infty m(k_1^*,s)ds,\eqno(26)$$
$$ h_1(x)\le \|f\|_1^{}\sup_y k_r(x,y)\le \alpha z r,\eqno(27)$$
so that letting $r=k_1^*(z)$ in (26) and (27) leads to
$$\eqalign{(Tf)^{**}(t)&={1\over t}\int_0^t (Tf)^*\le \|Tf\|_\infty^{}\le \|h_1\|_\infty^{}+ \|h_2\|_\infty^{}\cr&\le
\alpha z \,k_1^*(z)+\alpha\int_{k_1^*(z)}^\infty m(k_1^*,s)ds=\alpha\int_0^z k_1^*(s)ds=\alpha z \,k_1^{**}(z),\cr}$$
which is (23). 
If in addition $k_2^*$ is integrable, then 
$$\eqalign{&\int_N h_2(x)d\nu(x)=\int_N d\nu(x)\im k^r(x,y)f(y)d\mu(y)\cr&=\im f(y)\bigg(\int_N k^r(x,y)d\nu(x)\bigg)d\mu(y)
\le \|f\|_1^{}\int_r^\infty m(k_2^*,s)ds\le \alpha z \int_r^\infty m(k_2^*,s)ds,\cr}\eqno(28)$$
therefore,   letting $r=k_2^*(t)$ and using (27) and (28)
$$\eqalign{t\,(Tf)^{**}(t)&\le \int_0^t h_1^*+\int_0^t h_2^* \le t\,\|h_1\|_\infty^{} +\int_0^\infty h_2^*\cr& \le t\,\alpha z\, k_2^*(t)+\alpha z\int_{k_2^*(t)}^\infty m(k_2^*,s)ds=\alpha z\,t\, k_2^{**}(t)\cr}$$
and this concludes the proof of  (23) and (25), in case both $k_1^*$ and $k_2^*$ are integrable.
If conditions (17) are assumed instead, then (23) still holds (since only integrability of $k_1^*$ was needed), and  
estimate (24) is an immediate consequence of the weak-type estimate (22).
 \endpf

\bigskip\noin{\bf Remark.} We emphasize here the new elements appearing in the lemma, as compared to O'Neil's original version. First, the  role of the two measures, as reflected in the explicit dependence on $k_1^*$ and $k_2^*$, and their bounds. Secondly,  the fact that O'Neil's lemma is really a two-variable statement; this is hinted in the Claim, even in O'Neil's original version,  but it does not seem to have been noticed before. Our original version of the lemma was just (20) with  $\tau=t$ which was suitable to prove Theorem 1 when $\beta_0=\beta$ (our first version) but not for $\beta_0<\beta$. The further improvements of O'Neil's lemma came about in our attempts to incorporate some of Cianchi's main results [Ci1] in our general framework (see Theorem 6).

\vskip1em

\pf Proof of Theorem 1. It is enough to assume that $k$ is  nonnegative, and show that for each nonnegative $f\in L^{\beta'}(M)$    with $\,\|f\|_{\beta'}^{}\le 1$ we have 
$$\int_N \exp\bigg[{\beta_0\over A\beta}(Tf)^\beta\,\bigg]d\nu\le C\eqno(29)$$
for some $C$ independent of $f$.

Pick any $p$ as in (18). By (19) of the improved   O'Neil's Lemma 2, with $\tau=t^{\beta/\beta_0}$ 
$$\eqalign{&(Tf)^*(t)\le(Tf)^{**}(t)\le C t^{-{1\over q}}\int_0^{t^{\beta/\beta_0}}f^*(u) u^{-1+{1\over p}}du+\int_{t^{\beta/\beta_0}}^\infty k_1^*(u) f^*(u)du\cr &=
C t^{-{1\over q}}\int_0^{t}f^*\big(u^{\beta\over\beta_0}\big) u^{-1+{\beta\over p\beta_0}}du+{\beta\over\beta_0}\int_t^\infty k_1^*\big(u^{\beta\over\beta_0}\big) f^*\big(u^{\beta\over\beta_0}\big)u^{{\beta\over\beta_0}-1}du.\cr}\eqno(30)$$
By Fact 4, combined with the fact that $k_j^*(t)=0$ for $t\ge\max\{\nu(N),\mu(M)\}$,
$$k_1^*\big(u^{\beta\over\beta_0}\big) \le A^{1\over\beta}u^{-{1\over\beta_0}}\big(1+C(1+|\log u|)^{-\gamma}\big)\,,\quad u>0\eqno(31)$$
(C denotes a positive constant that may change from place to place).

Combining (30) and (31) yields
$$ \eqalign{(Tf)^{**}(t)\le C & t^{-{1\over q}}\int_0^{t}f^*\big(u^{\beta\over\beta_0}\big) u^{-1+{\beta\over p\beta_0}}du
+\cr&+{\beta\over\beta_0}\int_t^{\mu(M)^{\beta_0/\beta}}\!\!A^{1/\beta}\big(1+C(1+|\log u|)^{-\gamma}\big)f^*\big(u^{\beta\over\beta_0}\big)u^{{\beta\over\beta_0}-1}du\cr}
$$
and therefore, with $t_1=\max\{\nu(N),\mu(M)^{\beta_0/\beta}\}$,
$$\eqalign{&\int_N \exp\bigg[{\beta_0\over A\beta}(Tf)^\beta\,\bigg]d\nu(x)=\int_0^{\nu(N)} \exp\bigg[{\beta_0\over A\beta}\big((Tf)^*(t)\big)^\beta\bigg]dt
\le\int_0^{\nu(N)}\exp\bigg[{\beta_0\over A\beta}\Big((Tf)^{**}(t)\Big)^\beta\Big]dt\cr&\le\int_0^{t_1} \exp\bigg[\bigg( C t^{-{1\over q}}\int_0^{t}f^*\big(u^{\beta\over\beta_0}\big) u^{-1+{\beta\over p\beta_0}}du+\cr&\hskip3em      
+ \bigg({\beta\over\beta_0}\bigg)^{1\over\beta'}\int_t^{t_1}\big(1+C(1+|\log u|)^{-\gamma}\big)f^*\big(u^{\beta\over\beta_0}\big)u^{{\beta\over\beta_0}-1}du\bigg)^\beta\,\bigg]dt.\cr}$$

Now we make the changes of variables $u=e^{-x},\, t=e^{-y}$, and we let $y_1=-\log t_1$ and 
$$\phi(x)=\bigg({\beta\over\beta_0}\bigg)^{1\over\beta'}f^*\big(e^{-{\beta x\over\beta_0}}\big)\,e^{-{\beta-1\over\beta_0}x}.$$ Notice  that $\phi$ is defined  on $[y_1,\infty)$ and $\|\phi\|_{\beta'}^{}=\|f^*\|_{\beta'}^{}=\|f\|_{\beta'}^{}\le1$.

With these changes, estimate (29) reduces to
$$\int_{y_1}^\infty \exp\bigg[\bigg( H\int_{y}^\infty \phi(x) e^{y-x\over q}dx+\int_{y_1}^y\Big(1+H(1+|x|)^{-\gamma}\Big)\phi(x)dx\bigg)^\beta-y\bigg]\,dy\le C\eqno(32)$$
where $H$ is a suitable, fixed, positive constant. 

Define 
$$g(x,y)=\cases{1+H(1+|x|)^{-\gamma} & if $y_1\le x\le y$ \cr \cr
 H e^{{y-x\over q}} & if $y_1\le y<x<\infty$ \cr}\eqno(33)$$
and  
$$F(y)=y-\bigg(\int_{y_1}^\infty g(x,y)\phi(x) dx\bigg)^\beta\eqno(34)$$
which is defined for $y\in[y_1,\infty)$. Estimate (32) is a direct consequence of the following modified Adams-Garsia's lemma:

\bigskip
\proclaim Lemma 3. Suppose that $\phi:[y_1,\infty)\to [0,\infty)$ satisfies $\ds{\int_{y_1}^\infty \phi^{\beta'}\le 1},$ and  $g$ and $F$ be defined as in (33), (34), with $H>0$, $\beta>1$, $q>0$ 
 and $\ds{{1\over\beta}+{1\over\beta'}=1}$. Then there exists 
a constant $C$ independent   of $\phi$ such that 
$$\int_{y_1}^\infty e^{-F(y)}dy\le C.\eqno(35)$$
\par

This lemma  differs from the original Adams-Garsia
 lemma (Lemma 1 in [Ad1]) by the perturbation term $H(1+|x|)^{-\gamma}$ for $x\le y$ (which  was not present in  Adams-Garsia's lemma). In his original work Moser had 1 for $x\le y$ and 0 for $x>y$ which makes the argument much simpler. The proof below is a modification of the proof or Lemma 3.2 in [F], which was itself a modification of the proof of Lemma 1 in [Ad1]. We note that in [FFV] there is an even more general version of Lemma 3, which appeared after that in [F], but we decided to include its  proof in order to make our results  self-contained.

\bigskip
\pf Proof of Lemma 3. Let $E_\lambda=\{y\ge y_1:\, F(y)\le \lambda\}$ and let $|E_\lambda|$ be its Lebesgue measure.  Then
$$\int_{y_1}^\infty e^{-F(y)}dy=\int_{-\infty}^\infty |E_\lambda| e^{-\lambda}d\lambda.$$

\proclaim Claim 1. There exists $c\ge0$ independent of $\phi$ such that if $E_\lambda\neq \emptyset$, then $\lambda\ge -c$, i.e.
$\inf_{y\ge y_1} F(y)\ge -c>-\infty$.\par \medskip
\proclaim Claim 2. There exists $C$ independent of $\phi$ and $\lambda$ such that for every $\lambda\in \R$

$$|E_\lambda|\le C(1+|\lambda|).\eqno(36)$$
\par

Claims 1 and 2 imply (35) since 
$$\int_{y_1}^\infty e^{-F(y)}dy=\int_{-c}^\infty |E_\lambda| e^{-\lambda}d\lambda\le C\int_{-c}^\infty (1+|\lambda|)e^{-\lambda}d\lambda,$$
which is a constant independent of $\phi$.

\medskip
\pf Proof of Claim 1.  It is enough to assume that $\lambda<0$ and $y_1-\lambda>0$. If $y\in E_\lambda$ then  
$$\eqalign{&(y-\lambda)^{1\over\beta}\le\int_{y_1}^y\Big(1+H(1+|x|)^{-\gamma}\Big)\phi(x)dx+H\int_y^\infty e^{y-x\over q}\phi(x)dx\cr
& \le \bigg(\int_{y_1}^y \phi^{\beta'}\bigg)^{1\over\beta'}\bigg(\int_{y_1}^y \Big(1+H(1+|x|)^{-\gamma}\Big)^\beta dx\bigg)^{{1\over\beta}}+H
\bigg(\int_y^\infty\phi^{\beta'}\bigg)^{1\over\beta'}\bigg(\int_y^\infty e^{(y-x){\beta\over q}}dx\bigg)^{1\over\beta}.\cr}$$
Note that for $a,b\ge0$ and $\beta\ge1$ 
$$(a+b)^\beta\le a^{\beta}+\beta 2^{\beta-1}(a^{\beta-1}b+b^{\beta})\eqno(37)$$
(identity at $b=0$, and $b-$derivative of LHS smaller than $b-$derivative of RHS). Hence
$$\int_{y_1}^y \Big(1+H(1+|x|)^{-\gamma}\Big)^\beta dx\le\int_{y_1}^y \Big(1+H_1(1+|x|)^{-\gamma}\Big) dx\le y-y_1+d_1= y+d$$
some $d\in\R$, independent of $y$ (here is where we use $\gamma>1$).

As a result, if we let 
$$L(y)=\bigg(\int_y^\infty \phi^{\beta'}\bigg)^{1\over\beta'}\in [0,1].$$
 we have (using (37) again)
$$\eqalign{y-\lambda&\le \Big[\big(1-L(y)^{\beta'}\big)^{{1\over\beta'}}(y+d)^{{1\over\beta}}+CL(y)\Big]^\beta\cr&\le
\big(1-L(y)^{\beta'}\big)^{{\beta\over\beta'}}(y+d)+\beta 2^{\beta-1}\Big[\big(1-L(y)^{\beta'}\big)^{\beta-1\over\beta'}(y+d)^{\beta-1\over\beta} CL(y)+C^\beta L(y)^\beta\Big]\cr}$$

Since $\beta,\beta'>1, \,L(y)\in[0,1]$ and   $\big(1-L(y)^{\beta'}\big)^{\beta\over\beta'}\le 1-{\ds{1\over\beta'}}L(y)^{\beta'}$,
 if we let   $z=(y+d)^{1/\beta'}L(y)\ge0$ we obtain
$$z^{\beta'}\le D z+\beta' \lambda+D$$
for some constant $D$ (independent of $y$ and $\phi$).
Since $z^{\beta'}-Dz-D$ has a finite negative minimum on $[0,\infty)$, we deduce that if $E_\lambda\neq\emptyset$ then 
$\lambda\ge -c$, for some $c\ge0$ (independent of $y$ and $\phi$). 

 Note also that for large $z$ we have $Dz\le {1\over2} z^{\beta'}$ so that  $z^{\beta'}\le C(|\lambda|+1)$
or
$$ (y+d)^{{1\over\beta'}}L(y)\le C(|\lambda|^{{1\over\beta'}}+1)\eqno(38)$$
for some $C$ independent of $y$, $\phi$, and $\lambda$.
\bigskip
\pf Proof of Claim 2. 
It is enough to  prove that there exists $H>0$ (independent of $\phi$) such that for any $\lambda\in\R$
$$t_1,t_2\in E_\lambda\, {\hbox { and }}\, t_2>t_1>H|\lambda|+H\,\Longrightarrow \, t_2-t_1\le H|\lambda|+H.\eqno(39)$$
Indeed, if this is the case, then (recall that $E_\lambda\subseteq[y_1,\infty)$)
$$\eqalign{|E_\lambda|&=\big|E_\lambda\cap\{t: t\le H|\lambda|+H\}\big|+\big|E_\lambda\cap\{t:t>H|\lambda|+H\}\big|\cr& \le H|\lambda|+H-y_1+\sup_{t_2>t_1>H|\lambda|+H\atop t_1,t_2\in E_\lambda} (t_2-t_1)\le C|\lambda|+C.\cr} $$

To show (39), pick $t_1,t_2\in E_\lambda$, $\;t_2>t_1$, so that, arguing  as in the proof of Claim 1
$$\eqalign{(t_2-\lambda)^{{1\over\beta}}&\le \int_{y_1}^\infty g(x,t_2)\phi(x)dx\le \bigg(\int_{y_1}^{t_1}g(x,t_2)^\beta\bigg)^{{1\over\beta}}\bigg(\int_{y_1}^{t_1}\phi^{\beta'}\bigg)^{{1\over\beta'}}\cr& \hskip2em + \bigg(\int_{t_1}^{t_2}g(x,t_2)^\beta\bigg)^{{1\over\beta}}\bigg(\int_{t_1}^{t_2}\phi^{\beta'}\bigg)^{{1\over\beta'}}+\bigg(\int_{t_2}^{\infty}g(x,t_2)^\beta\bigg)^{1\over\beta}\bigg(\int_{t_2}^{\infty}\phi^{\beta'}\bigg)^{{1\over\beta'}}\cr&
\le (t_1+d)^{{1\over\beta}}+(t_2-t_1+d_1)^{{1\over\beta}}\bigg(\int_{t_1}^\infty\phi^{\beta'}\bigg)^{{1\over\beta'}}+C\bigg(\int_{t_1}^\infty\phi^{\beta'}\bigg)^{{1\over\beta'}}\cr &=(t_1+d)^{{1\over\beta}}+\big((t_2-t_1+d)^{{1\over\beta}}+C\big) L(t_1)\cr}$$
which , using (37) and (38), implies
$$\eqalign{t_2-\lambda&\le t_1+d+\beta 2^{\beta-1}\bigg[(t_1+d)^{\beta-1\over\beta}\big((t_2-t_1)^{{1\over\beta}}+C\big)L(t_1)+\big((t_2-t_1)^{{1\over\beta}}+C\big)^\beta L(t_1)^\beta\bigg]\cr &\le
 t_1+d+\beta 2^{\beta-1}\bigg[\big((t_2-t_1)^{{1\over\beta}}+C\big)(t_1+d)^{{1\over\beta'}}L(t_1)+2^\beta(t_2-t_1)L(t_1)^\beta+2^\beta C^\beta\bigg]\cr&\le
 t_1+\big((t_2-t_1)^{{1\over\beta}}+C\big)(C|\lambda|^{{1\over\beta'}}+C)+C(t_2-t_1)L(t_1)^\beta+C\cr&
\le t_1+{t_2-t_1\over\beta}+ {(C|\lambda|^{{1\over\beta'}}+C)^{\beta'}\over\beta'}+C(t_2-t_1)L(t_1)^\beta+C|\lambda|^{{1\over\beta'}}+C.\cr}$$

Hence, 
$$ {t_2-t_1\over\beta'}\le C|\lambda|+C+C(t_2-t_1)L(t_1)^\beta\le C|\lambda|+C+(t_2-t_1){C|\lambda|+C\over t_1+d}$$
so it follows that there is $C$ so that 
$$t_1+d>2\beta' C|\lambda|+2\beta' C\,\Longrightarrow \, t_2-t_1\le 2\beta'C|\lambda|+2\beta' C,$$ which is (39). 
Claim 2, Lemma 3 and Theorem 1 are thus completely proven.\endpf

\bigskip\medskip
\noin{\bf 2. Conditions for sharpness}
\bigskip
In the following theorem we prove that,   under suitable ``geometric" conditions, equality in (13),  implies that $\ds{\beta_0\over A\beta}$
in (15) or (16) is sharp, i.e. it cannot be replaced by a larger constant. We state and prove the general vector-valued case, since it does not follow directly from the scalar case, as opposed to the proof of Theorem 1'. It will be apparent from the proof that the same result will also hold for complex-valued operators (see Remark 1 after the proof of Theorem 4).
\smallskip
For measurable $F:M\to\R^n$ and  $K:N\times M\to\bar\R^n$ let  
$$TF(y)=\int_M K(x,y)\cdot F(y)\,d\mu(y)$$
if the integral is well defined.

\bigskip\eject

\proclaim Theorem 4. Suppose that $k(x,y)=|K(x,y)|$ satisfies
$$ m(k_1^*,s)=A s^{-\beta}\big(1+O(\log^{-\gamma} s)\big),\qquad {\hbox { as }} s\to+\infty,\eqno(40)$$
or equivalently 
$$ k_1^*(t)=A^{1/\beta} t^{-1/\beta} \Big(1+O\big(|\log t|^{-\gamma}\big)\Big),\qquad {\hbox { as }} t\to0,\eqno(41)$$
 and 
$$ m(k_2^*,s)\le B s^{-\beta_0}$$
as $s\to+\infty$,  for some $\beta,\gamma>1$, $\,0<\beta_0\le \beta$ and $B>0$.
Suppose that there exist $x_m\in N$, measurable sets $B_m\subseteq N,\, E_m\subseteq M$, $\,m\in \N$, with the following properties:\medskip
\noin {a)} $E_m\supseteq\{y:\,|K(x_m,y)|>m\},$  $\;\mu(E_m)=O(m^{-\beta})$, as $m\to \infty$\smallskip
\noin{b)} there exist constants $c_1,c_2>0$ such that $c_1m^{-\beta_0}\le \nu(B_m)\le  c_2m^{-\beta_0},\, m=1,2....$
\smallskip
\noin{c)} 
$$k^*(x_m,\cdot)(t)\ge A^{1/\beta}t^{-1/\beta}\Big(1-c_3|\log t|^{-\gamma}\Big),\quad 0<t<t_0<1\eqno(42)$$
\noin{d)}
$$ \int_{M\setminus E_m} \big|\big(K(x,y)- K(x_m,y)\big)\cdot K(x_m,y)\big|\, |K(x_m,y)|^{\beta-2} d\mu(y)\le c_4\,,\quad \forall x\in B_m
\eqno(43)$$
with $c_3,c_4$ independent of $x$ and $m$. Then, the Adams inequality (16) holds and it is sharp, in the sense that 
$$\sup_{F\in L^{\beta'}(M)} \int_N \exp \bigg[\alpha\,\bigg({|TF|\over\|F\|_{\beta'}^{}}\bigg)^\beta\,\bigg]d\nu=+\infty,\qquad \forall \alpha>{\beta_0\over A\beta}.$$
More specifically, if a), b) , c) hold and 
$$\Phi_m(y)=K(x_m,y)|K(x_m,y)|^{\beta-2}\chi_{M\setminus E_m}^{}(y)\eqno(44)$$ then $\Phi_m\in
 L^{\beta'}$ with 
$$\|\Phi_m\|_{\beta'}^{\beta'}=A\,\log{1\over\mu(E_m)}+O(1),\eqno(45)$$
and if d) also holds then
$$\lim_{m\to\infty}\int_N \exp \bigg[\alpha\,\bigg({|T\Phi_m|\over\|\Phi_m\|_{\beta'}^{}}\bigg)^\beta\,\bigg]d\nu=+\infty,
\qquad \forall \alpha>{\beta_0\over A\beta}.\eqno(46) $$
\par

\bigskip\eject
\noin{\bf Remarks.} \smallskip

\noindent{\bf 1.} If there is a point $x_0$ such that $k_1^*(t)=k^*(x_0,\cdot)(t)$ for small $t$, then typically one can choose $x_m=x_0$, so that (42) is automatically true. In the context of metric spaces one can typically choose $E_m$ to be the $m-$th level set of $|K(x_0,y)|$, or a possibly slightly larger set, and $B_m$  a suitable small ball around $x_0$.  In all the applications we know, the only minor technical check  is about  the integral estimate in (43), which  is usually a consequence of H\"older continuity  estimates on $K(x,y)$. This point is illustrated clearly in all the applications presented in  Section 5.

\smallskip
\noin{\bf 2.} In the scalar case  $K(x,y)=k(x,y)$ condition d) obviously becomes
$$\int_{M\setminus E_m} |k(x,y)- k(x_m,y)|\, |k(x_m,y)|^{\beta-1} d\mu(y)\le c_4\,,\qquad \forall x\in B_m.\eqno(47)$$
In the vector-valued case condition d) is implied by 
$$\int_{M\setminus E_m}|K(x,y)- K(x_m,y)|\, |K(x_m,y)|^{\beta-1} d\mu(y)\le c_4\,,\qquad \forall x\in B_m.$$
\smallskip
\noin{\bf 3.} The classical form of a  Moser-Trudinger inequality for a differential (or pseudodifferential) operator of order $d$ takes the form
$$\int_N \exp \bigg[\alpha\,\bigg({|u|\over\|Pu\|_{p}^{}}\bigg)^{p'}\,\bigg]d\nu\le C\eqno(48)$$
where $P$  acts on a suitable subspace of $L^p(N)$ (usually a Sobolev space). A  lower bound for $\alpha$ can be achieved via a 
 representation formula $u=T(Pu)$, where $T$ is an integral operator with kernel $K$, satisfying the hypothesis of Theorem 1 or 1'. If  the conditions in Theorem 4 are satisfied, then the sharpness of the constant follows immediately if one is able to produce a sequence $u_m$ in the given space such that $Pu_m=\Phi_m$, the extremizing sequence of Theorem 4. When dealing with scalar functions this is usually possible (see Theorems 6 and 10). Another similar  way to obtain an upper bound for $\alpha$ is  by choosing a suitable sequence of functions $u_m$ and sets $B_m\subseteq N$ such that $u_m\ge \delta_m$ on $B_m$, via the inequality 
$$\alpha\le \liminf_n \bigg({\|Pu_m\|_p\over\delta_m}\bigg)^{p'}\log{1\over \nu(B_m)}\eqno(49)$$
which follows easily from (48). This approach is slightly more flexible in that the  $u_m$ may not be the exact inverse images of the $\Phi_m$,  even though they usually differ from those by  negligible terms. 


\bigskip
The following lemma will play an important  role in the proof of Theorem 4:
\smallskip
\proclaim Lemma 5. Let $f:M\to \R$ be measurable, and $E\subseteq M$ measurable with $0<\mu(E)<\mu(M)$. Let 
$$\wtilde f(y)=\cases{\ds\mathop\essup\limits_{z\in M\setminus E} |f(z)|& if $\,y\in E$\cr
 f(y) & if $\,y\in M\setminus E$.\cr}$$
Then
$$\wtilde f^*(t)\ge f^*(t),\qquad \mu(E)\le t\le \mu(M).$$
Moreover,
$$\int_{M\setminus E} |\wtilde f|^\beta=\int_{\mu(E)}^{\mu(M)} [\wtilde f^*(t)]^\beta dt.$$
\par\medskip
\pf Proof. Suppose first that $f$ is essentially bounded on $M\setminus E$ (actually this is all we need for the proof of Theorem 4). If $s_0=\ds\mathop\essup\limits_{z\in M\setminus E} |f(z)|$, then 
$|f|\le s_0$ a.e. on $M\setminus E$, so that  $m(f,s_0)\le \mu(E)$. This implies
that for $ \mu(E)\le t\le  \mu(M)$ $$f^*(t)\le f^*\big(\mu(E)\big)=\inf\{s:m(f,s)\le \mu(E)\}\le s_0$$
which proves the claim if $\wtilde f^*(t)=s_0$ (it cannot be $>s_0)$. On the other hand
$$\{y:\,|\wtilde f(y)|>s\}=\{y:\,|f(y)|>s\}\cup E,\,\qquad 0<s<s_0$$
 hence $m(\wtilde f,s)\ge m(f,s)$ if $0<s<s_0$, so that if $\wtilde f^*(t)=\inf\{s: m(\wtilde f,s)\le t\}<s_0$, then 
$$\wtilde f^*(t)=\inf\{s<s_0: m(\wtilde f,s)\le t\}\ge\inf\{s<s_0: m(f,s)\le t\}\ge f^*(t).$$

For the last statement, note that $\wtilde f^*(t)=s_0$ for $0<t<\mu(E)$ (indeed $m(\wtilde f,s)=0$ for $s\ge s_0$, and for any $s<s_0$ we have $|\wtilde f(y)|>s$ on $E$  and on a set of positive measure inside $M\setminus E$, i.e. $m(\wtilde f,s)>\mu(E)$, for $s<s_0$). Hence
$$\int_M  |\wtilde f|^\beta=\int_{\mu(E)}^{\mu(M)} [\wtilde f^*(t)]^\beta dt+s_0^\beta\mu(E)$$
but also 
$$\int_M |\wtilde f|^\beta =\int_{M\setminus E}|\wtilde f(y)|^\beta d\mu(y)+s_0^\beta \mu(E).$$
A standard approximation argument (truncation and the monotone convergence theorem) completes the proof for general $f$.
\endpf
\pf Proof of Theorem 4.  Define for $m\in \N$ 
$$\Phi_m(y)=K(x_m,y)|K(x_m,y)|^{\beta/\beta'-1}\chi_{M\setminus E_m}^{}(y)$$
$$F_m=\{y: \,|K(x_m,y)|>m\}\subseteq E_m$$
so that, by (41) and a), and since $F_m$ is a level set for $|K(x_m,y)|$, 
$$\eqalign{\|\Phi_m\|_{\beta'}^{\beta'}&=\int_{M\setminus E_m}|K(x_m,y)|^{\beta}d\mu(y)\le \int_{M\setminus F_m}|K(x_m,y)|^{\beta}d\mu(y)=
\int_{\mu(F_m)}^{\mu(M)} [k^*(x_m,\cdot)(t)]^\beta dt\cr&\le \int_{\mu(F_m)}^{\mu(M)} A \Big(1+C(1+|\log t|)^{-\gamma}\Big)\,{dt\over t}=A\,\log{1\over\mu(F_m)}+C\le A\,\log{1\over\mu(E_m)}+C'\cr}$$
(the last inequality follows from the assumptions a) and  c)).
On the other hand, if we define 
 $$\km(y)=\cases{{\ds{\rm ess} \!\!{\sup_{\!\!\!\!\!\!z\in M\setminus E_m} }}|K(x_m,z)| & if $\,y\in E_m$\cr\cr |K(x_m,y)| &if $\,y\in M\setminus E_m$\cr}$$
then by Lemma 5 we have $\km^*(t)\ge k^*(x_m,\cdot)(t)$, for $\mu(E_m)\le t\le \mu(M)$, so that (by (42))
$$\eqalignno{\|\Phi_m\|_{\beta'}^{\beta'}&=\int_{M\setminus E_m}\!\!|K(x_m,y)|^{\beta}d\mu(y)= \int_{M\setminus E_m}\!\!|\km(y)|^{\beta}d\mu(y)
=\int_{\mu(E_m)}^{\mu(M)} [\km^*(x_m,\cdot)(t)]^\beta dt\cr&\ge \int_{\mu(E_m)}^{\mu(M)} A \Big(1-C(1+|\log t|)^{-\gamma}\Big)\,{dt\over t}= A\,\log{1\over\mu(E_m)}-C.&(50)\cr} $$
which gives (45).
Now, for $x\in B_m$, using (50) and (43)
$$\eqalign{&T\Phi_m(x)=\im K(x,y)\cdot\Phi_m(y)\,d\mu(y)=\int_{M\setminus E_m} K(x,y)\cdot K(x_m,y)\,|K(x_m,y)|^{\beta/\beta'-1}
d\mu(y)\cr&=\int_{M\setminus E_m} |K(x_m,y)|^{1+\beta/\beta'}d\mu(y)+\int_{M\setminus E_m}\Big(K(x,y)- K(x_m,y)\Big)\cdot K(x_m,y)|K(x_m,y)|^{\beta/\beta'-1}\,d\mu(y)\cr& \ge A\log{1\over \mu(E_m)}- C
\cr}$$
with $C$ independent of $m$. Hence, if $\,\wtilde \Phi_m=\Phi_m  \|\Phi_m\|_{\beta'}^{-1}$
and  $x\in B_m$
$$ T\wtilde \Phi_m(x)\ge {A\log\ds{1\over \mu(E_m)}- C\over \Big(A\,\log\ds{1\over\mu(E_m)}\Big)^{1/\beta'}+O(1)}=\Big(A\log{1\over \mu(E_m)}\Big)^{1/\beta}+O(1).$$
Finally, if $\alpha> \ds{\beta_0\over A\beta}$
$$\eqalign{\int_N \exp\Big[\alpha |T\wtilde \Phi_m(x)|^\beta\Big]d\nu(x)&\ge \int_{B_m}e^c \exp\Big[\alpha A\log{1\over \mu(E_m)}\Big] \,d\nu(x)\cr& =e^c\nu(B_m)\big(\mu(E_m)\big)^{-\alpha A}\ge C m^{-\beta_0+\alpha A \beta}\to+\infty.\cr}$$
\endpf
\bigskip
\noin{\bf Remark 1.}  It is clear from the proof just shown that Theorem 4 holds almost verbatim when $K$ is complex-valued and $T$  acts on complex-valued functions. The functions $\Phi_m$  need only to be replaced by
$$\Phi_m(y)=\bar {K(x_m,y)}\,|K(x_m,y)|^{\beta -2}\chi_{M\setminus E_m}^{}(y)$$
\vskip1em
\noin{\bf 3. Sharpness in $\gamma$}
\bigskip
In this section we show that if $\gamma\le 1$ in (13) then  the conclusion of Theorem 1 is in general false. We do this by considering the simplest setting,  namely $N=M=B(0,1)=\{x\in\Rn:\,|x|\le1\}$, $0<d<n,\,\beta'={\ds {n\over d}},\, \beta=\beta_0={\ds{n\over n- d}}$, and for the operator $T_d$ defined as
$$T_ d f(x)=\int_{B(0,1)} |x-y|^{ d-n}\bigg(1+{1\over 1+\big|\log|x-y|\big|}\bigg) f(y)dy\,\qquad f\in L^{\beta'}\big(B(0,1)\big).$$
As we mentioned previously, in this case $A=\ds{\omega_{n-1}\over n}.$ For $0<r<1$ consider $f_r(x)=|x|^{- d}\chi_{\{r\le|x|\le1\}}^{}(x)$, so that $\|f_r\|_{\beta'}^{}=\big(\omega_{n-1}\log\ds{\ts {1\over r}}\big)^{1/\beta'}$. Letting $\wtilde f_r(x)=f_r(x)\|f_r\|_{\beta'}^{-1}$ we obtain
$$\eqalign{\big(\omega_{n-1}&\log\ds{\ts {1\over r}}\big)^{1/\beta'}T_ d\wtilde f_r(x)=\int_r^1s^{n-1- d} ds\int_{S^{n-1}} 
\big||x|e_1-s\sigma\big|^{ d-n}\bigg(1+{1\over1+\Big|\log\big||x|e_1-s\sigma\big|\Big|}\bigg)d\sigma \cr
&=\int_{|x|}^{|x|/r} {dt\over t}\int_{S^{n-1}} 
|te_1-\sigma|^{ d-n}\bigg(1+{1\over1+\big|\log|te_1-\sigma|+\log|x|-\log t\big|}\bigg)d\sigma
\cr
& =\int_{|x|}^{|x|/r} {dt\over t}\int_{S^{n-1}} 
\bigg(1+{1\over1+\big|\log|te_1-\sigma|+\log|x|-\log t\big|}\bigg)d\sigma+\cr&\hskip2em+ \int_{|x|}^{|x|/r} {dt\over t}\int_{S^{n-1}} 
\big(|te_1-\sigma|^{ d-n}-1\big)\bigg(1+{1\over1+\big|\log|te_1-\sigma|+\log|x|-\log t\big|}\bigg)d\sigma.\cr
}$$

\noin If $|x|\le r/2$ and  $|x|\le t\le |x|/r$, then $\log{1\over2}\le \log|te_1-\sigma|\le\log{3\over2}$ and $\big||te_1-\sigma|^{ d-n}-1\big|\le Ct,$ 
so that the second  term above is bounded, for $|x|\le r/2$, and 
$$\eqalign{T_ d \wtilde f_r(x)&\ge{\omega_{n-1}^{1/\beta}\over \big(\log{\ts {\ts {1\over r}}}\big)^{1/\beta'}}\bigg[\int_{|x|}^{|x|/r}\bigg(1+{1\over 1+\log2+\log t-\log|x|}\bigg){dt\over t}-C\bigg]\cr&\hskip-3em\ge {\omega_{n-1}^{1/\beta}\over \big(\log{\ts {1\over r}}\big)^{1/\beta'}}\bigg[\log{\ts {1\over r}}+\log\bigg(1+{\log{\ts {1\over r}}\over1+\log 2}\bigg)-C\bigg]\ge \big(\omega_{n-1}\log{\ts {1\over r}}\big)^{1/\beta}\bigg[1+{\log \log{\ts {1\over r}}\over2\log{\ts {1\over r}}}\bigg]
\cr}$$
for $0<r\le r_0<1$, with $r_0$ small enough. Finally, for $r\le r_0$ 
$$\eqalign{\int_{B(0,1)} \exp\bigg[{n\over \omega_{n-1}}&\big(T_ d \wtilde f_r\big)^\beta\bigg] dx\ge\int_{B(0,r/2)}
\exp\bigg[n\log{1\over r}\bigg(1+{\log \log{\ts {1\over r}}\over2\log{\ts {1\over r}}}\bigg)^\beta\bigg]dx\cr
&\ge {\omega_{n-1}\over n}\Big({r\over2}\Big)^n r^{-n} \exp\Big[{\ts{1\over2}}n\beta\log\log{1\over r}\Big]={\omega_{n-1}\over n}\, 2^{-n} \big(\log\ts{1\over r}\big)^{n\beta/2}\to+\infty,\cr
}$$
as $r\to0$.

\bigskip
\bigskip
\centerline{\bf PART II: APPLICATIONS}
\bigskip\bigskip
\noin{\bf 4. Adams and Moser-Trudinger trace inequalities}
\bigskip

As a first illustration of our theorems we give a simultaneous extension of Adams' original results in [Ad1], and Cianchi's recent sharp Moser-Trudinger trace inequality in [Ci1]
$$\int_\Omega \exp\bigg[\lambda\omega_{n-1}^{1\over n-1}\bigg({|u(x)|\over \|\nabla u \|_{L^n(\Omega)}}\bigg)^{n\over n-1}\bigg]\,d\nu(x)\le C\eqno(51)$$
valid for all $u\in W_0^n(\Omega)$, where $L^n(\Omega)$ and $W_0^n(\Omega)$ are w.r. to the Lebesgue measure, and where $\nu$ denotes a positive Borel measure on $\Omega$ such that 
$$\exists \lambda\in(0,n]\,\, {\rm and }\,\, r_0>0\,\,:\,\,\nu\Big(B(x,r)\cap \Omega\Big)\le C r^\lambda,\qquad \forall x\in \R^n,\,\, \forall r\in(0,r_0].\eqno(52)$$
Here and throughout the rest of this work 
$$B(x,r)=\{ y\in\R^n: |y-x|<r\}.$$
Note that (51) is itself an extension of Adams' original result when $\lambda=n$, in the special case of the gradient. Cianchi's results are given in the slightly more general framework of Lorentz-Sobolev spaces.

As an immediate consequence of our general theorem we will now derive  Cianchi's result, for arbitrary powers of the Laplacian and their gradients.   In other words, we will   extend Adams' result (2), in the case of measures satisfying (52).

In the following,  by  $L^p(\Omega)$ we mean the space of Lebesgue measurable functions $f$  on $\Omega$ such that $\|f\|_p:=\Big(\int_\Omega |f|^pdx\Big)^{1/p}<\infty$.
Also,  $\lap$ will always denote the positive
Laplacian on $\Rn$. The fractional powers of $\lap$ are defined on the Schwarz 
class $\S$ as
$$(\lap^{ d/2} \phi)\;{\what{\vphantom{S}}}\;(\xi)=(\;2\pi|\xi|\;)^ d\;\what\phi(\xi)\,,
\qquad\phi\in \S.$$
Here $ d\in\R$ and the Fourier transform is defined as  $\;\what\phi(\xi)=\irn 
e^{-2\pi i x\cdot\xi} \phi(x)dx$. By duality we can extend $\lap^{ d/2} $ to an operator
acting on $\S'$, the space of tempered distributions.

When $0< d<n$ and $f\in L^p(\Rn)$, $p\ge1$, the equation $\lap^{ d/2} \phi=f$ has a unique 
solution in $L^q$, with $q^{-1}=p^{-1}- d/n$. This solution is given
explicitly in terms of  the Riesz potential 
$$\lap^{- d/2}f(x)=c_ d\irn |x-y|^{ d-n}f(y)dy\eqno(53)$$
with
$$c_ d={\Gamma\big({n- d\over2}\big)\over 2^ d\pi^{n/2}\Gamma\big({ d\over2}\big)}.\eqno(54)$$ 

In other words, the distributional Fourier transform of the RHS
of (53) coincides with $(2\pi|\xi|)^{- d}\what f(\xi)$. 

The usual Sobolev space on $\R^n$ is denoted by $W^{ d,p}(\Rn)$ and the space $W_0^{ d,p}(\Omega)$ is the closure of $C_0^\infty(\Omega)$ in $W^{ d,p}(\Rn)$.
\medskip

\proclaim Theorem 6. Let $\Omega$ be open and bounded on $\Rn$, $n\ge3$, and let $\nu$ be a positive Borel  measure  on $\Omega$ satisfying (52). For $0< d<n$, $ d\in\N$, let $p=\ds{n\over d}$ and  $\,\ds{{1\over p}+{1\over p'}=1}$. Then, there exists $C$ such that
$$\int_\Omega \exp\bigg[{\lambda\,c_ d^{-p'}\over \omega_{n-1}}\bigg({|u(x)|\over\|\Delta^{d/2}u\|_p}\bigg)^{p'}\bigg]d\nu(x)\le C,\qquad  d {\hbox{ even}}\eqno(55)$$
and
$$\int_\Omega \exp\bigg[{\lambda\over\omega_{n-1}}\big(c_{ d+1}(n- d-1)\big)^{-p'}\bigg({|u(x)|\over\|\nabla\Delta^{{ d-1}\over2}u\|_p}\bigg)^{p'}\bigg]d\nu(x)\le C,\qquad  d {\hbox{ odd}}\eqno(56)$$
for all $u\in W_0^{ d,p}(\Omega)$. The constants appearing inside the exponents in (55) and (56) are sharp, provided there exists $x_0\in\Omega$ such that $\nu\big(B(x_0,r)\cap \Omega\big)\ge C_1 r^\lambda$ for $0< r\le r_1$, some $r_1,C_1>0$.
  \par

\bigskip

When $\lambda=n$ and $ d=m$ one recovers the constants $\beta(m,n)$ appearing in [Ad1]. When $ d=1$  the constant in 
(56) coincides with that of Cianchi, for $0<\lambda\le n$. 
\medskip
It is clear that it is enough to  prove the theorem if $u$ is smooth  with compact support inside $\Omega$. Secondly, for $ d$ even
$$u(x)=c_ d\int_{\Omega} |x-y|^{ d-n} \Delta^{ d/2}u(y)dy$$ 
and for $ d$ odd
$$u(x)=c_{ d+1}(n- d-1)\int_\Omega |x-y|^{ d-n-1}(x-y)\cdot\nabla \Delta^{{ d-1\over2}}u(y)dy,\eqno(57)$$ 
and therefore the inequalities of Theorem 6 are instant consequences of the following:

\proclaim Theorem 7. Let $\Omega$ be open and bounded on $\Rn$, $n\ge1$, and let $\nu$ be a positive Borel  measure  on $\Omega$ satisfying (52). For $0< d<n$, let $p=\ds{n\over d}$ and  $\,\ds{{1\over p}+{1\over p'}=1}$. Define for $f:\Omega\to\R$, $f\in L^p(\Omega)$
$$Tf(x)=\int_\Omega |x-y|^{ d-n} f(y)dy$$
and for $F:\Omega\to\R^n$, $F\in L^p(\Omega)$
$$T_1F(x)=\int_\Omega |x-y|^{ d-n-1}(x-y)\cdot F(y)dy.$$
Then, there exists $C$ such that
$$\int_\Omega \exp\bigg[{\lambda\over \omega_{n-1}}\bigg({|Tf(x)|\over\|f\|_p}\bigg)^{p'}\bigg]d\nu(x)\le C\eqno(58)$$
for every $f\in L^p(\Omega)$, and
$$\int_\Omega \exp\bigg[{\lambda\over\omega_{n-1}}\bigg({|T_1F(x)|\over\|F\|_p}\bigg)^{p'}\bigg]d\nu(x)\le C\eqno(59)$$
for every $F\in L^p(\Omega)$. The constant $\ds{\lambda\over \omega_{n-1}}$ in (58) and (59) is sharp, provided there exists $x_0\in\Omega$ such that $\nu\big(B(x_0,r)\cap \Omega\big)\ge C_1 r^\lambda$, for $0< r\le r_1$, some $r_1,C_1>0$.
 \par

\pf  Proof.
It is easy to check that if $k(x,y)=|x-y|^{ d-n}$ then for large $s$

$$m(k_1^*,s)={\omega_{n-1}\over n}s^{-{n\over n- d}}$$
and, using (52),
$$m(k_2^*,s)\le C s^{-{\lambda\over n- d}},$$
so that  Theorems 1 and 1' immediately imply (58), (59). To verify sharpness, according to Theorem 4, (and Remark 1 following it) first assume WLOG that $x_0=0\in \Omega$,  then take  $x_m=0\in\Omega,$ and   $m,R$ large enough so that 

$$ \{y\in\Omega:\,|K(0,y)|>m\,\}= B(0, m^{-p'/n})\subseteq\Omega\subseteq B(0,R),\eqno(60)$$
and let $$r_m=m^{-p'/n},\;\; E_m=B(0,r_m)\;\;B_m=B(0,\ts{1\over2}r_m)\eqno(61)$$
with either $K(x,y)= |x-y|^{ d-n} $ or $K(x,y)= |x-y|^{ d-n-1}(x-y)$. Conditions a), b), c) of Theorem 4 are met, with $\beta=n/(n- d)$ and $\beta_0=\lambda/(n- d)$, so all we need to check is d), i.e.

$$\int_{\Omega\setminus E_m} \big|\big(K(x,y)- K(0,y)\big)\cdot K(0,y)\big|\, |K(0,y)|^{p'-2} dy\le C\,,\qquad  |x|\le {r_m\over2}$$
for either kernel. If $K(x,y)= |x-y|^{ d-n} $ we need to check

$$\sup_{|x|\le r_m/2}\;\int_{r_m\le |y|\le R} |y|^{- d}\Big||x-y|^{ d-n}-|y|^{ d-n}\Big|dy\le C\eqno(62)$$
for some $C$ independent of $m$, but this  estimate is an   immediate consequence of 
$$|x-y|^{ d-n}\le |y|^{d-n}\Big|{x\over|y|}-{y\over|y|}\Big|^{ d-n}\le 2^{n- d} |y|^{ d-n}\eqno(63)$$
and
$$ |y|^{- d}\Big||x-y|^{ d-n}-|y|^{ d-n}\Big|\le C|y|^{-n}\bigg|1-\Big|{x\over|y|}-{y\over|y|}\Big|^{n- d}\bigg|\le C|x||y|^{-n-1},\eqno(64)$$
both valid for any $ d<n$ (even negative) and $|x|\le{r_m\over2},\;|y|\ge r_m$. 
\medskip
If instead $K(x,y)= |x-y|^{ d-n-1}(x-y)$ then we are reduced to
$$\sup_{|x|\le r_m/2}\;\int_{r_m\le |y|\le R} \Big||x-y|^{ d-n-1}(x\cdot y-|y|^2)+|y|^{d-n+1}\Big|\,|y|^{- d-1}dy\le C$$
which is implied by
$$\sup_{|x|\le r_m/2}\;\int_{r_m\le |y|\le R} |x|\,|y|^{- d}\,|x-y|^{d-n-1}dy\le C$$
$$\sup_{|x|\le r_m/2}\;\int_{r_m\le |y|\le R} \Big||x-y|^{ d-1-n}-|y|^{d-1-n}\Big|\,|y|^{-( d-1)}dy\le C,$$
and both of these are also easy consequences of (63) and (64). 

\endpf

\noin{\bf Proof of Theorem 6}. The inequalities of Theorem 6 follows from the formulas
$$u(x)=Tf(x)=\int_{\Omega} |x-y|^{d-n} f(y) dy,\qquad f=c_d\Delta^{d/2}u$$
for $d$ even, and 
$$u(x)=T_1F(x)=\int_\Omega |x-y|^{d-n-1}(x-y)\cdot F(y)dy,\qquad F(y)=c_{d+1}(n-d-1)\nabla\Delta^{d-1\over2}u$$
 for $d$ odd.

As far as the sharpness statement, the proof we give below is a slight modification of  Adams' original method. We include some details here since similar formulas  will be used later in the proof  of  Theorem 12. It is possible however to give a proof based on our Theorem 4 by using representation formulas via the Green function of $\Delta^{d/2}$ on the unit ball, with zero boundary conditions. Later in Theorem 10 we will give a  sharp inequality for more general (scalar) operators, and the proof of sharpness given there is a more direct application of Theorem 4, and applies also to the case $d$ even of the present theorem (see Remark 2 after the proof of Theorem 10).

Let $r_m\to0^+$ and $B_m=B(0,r_m)\subseteq \Omega$, for $m$ large enough.

 Next, pick any $\delta\ge0  $ and  $\varphi\in C^\infty([\delta,1+\delta])$ such that  $\varphi^{(k)}(\delta)=0$ for $0\le k\le \ell$, some $\ell\ge d$, and $\varphi(1+\delta)=1+\delta, \varphi'(1+\delta)=1$, $\varphi^{(k)}(1+\delta)=0$ for $2\le k\le \ell$.

Define for $m$ large enough 
$$u_m(y)=\cases{0 &  for $\;|y|\ge e^{-\delta}$\cr \varphi\Big(\log\ds{1\over |y|}\Big) & for $\;e^{-1-\delta}\le |y|<e^{-\delta}$\cr  \log\ds{1\over|y|} & for $\;e^{1+\delta} r_m\le |y|< e^{-1-\delta}$\cr   \log\ds{1\over r_m}-\varphi\Big({\log\ds{|y|\over r_m}\Big)} & for $\;e^\delta r_m\le|y|<e^{1+\delta} r_m$\cr 
\log\ds{1\over r_m}& for $\;|y|<e^\delta r_m.$\cr}\eqno(65)$$ 

 Then, if we pick $\delta$ large enough we have $u_m\in W_0^{d,p}(\Omega)$ and it is easy to check that for $d$ even
$$\Delta^{d/2}u_m(y)=\cases{ \ds{|y|^{-d}\over \omega_{n-1}c_d}& if 
$e^{1+\delta} r_m<|y|< e^{-1-\delta}$\cr\cr O(1)|y|^{-d} & if $e^\delta r_m<|y|<e^{1+\delta} r_m$\cr\cr O(1) & otherwise\cr}$$
and
$$\|\Delta^{d/2}u_m\|_p^{p'}={1\over\omega_{n-1} c_d^{p'}}\Big(\log{1\over r_m}\Big)^{p'/p}+O(1).$$
This means that if  $$\int_\Omega \exp\bigg[\alpha \,\bigg({|u(x)|\over\|\Delta^{d/2}u\|_p}\bigg)^{p'}\bigg]d\nu(x)\le C$$
holds for all $u\in W^{d,p}(\Omega)$, then according to (49)
$$\alpha\le \liminf_n\Bigg({\|\Delta^{d/2} u_m\|_p\over\log(1/r_m)}\bigg)^{p'}\log{1\over\nu(B_m)}={\lambda\over \omega_{n-1} c_d^{p'}}.$$
The proof in the case $d$ odd is completely similar, and based on the identity 
$$\nabla\Delta^{d-1\over2}\log{1\over|y|}={(d-1)y|y|^{-d-1}\over \omega_{n-1}c_{d-1}}=-{y|y|^{-d-1}\over \omega_{n-1}(n-d-1)c_{d+1}}.$$
\endpf

We remark here that the above theorem can be formulated and proved  in the setting of compact Riemannian manifolds, in the same spirit as in [F]. The point is that such theorem holds for any pseudodifferential operator of order $d$, whose leading symbol is $\Delta^{ d/2}$ (see next section for even a more a general result). For such operators the fundamental solution $k(P,Q)$ satisfies locally
$$k(P,Q)=c_ d d(P,Q)^{ d-n}\Big(1+O\big(d(P,Q)^{\e}\big)\Big)$$
for some $\e>0$, where $P,Q$ are points on the manifold, and $d(P,Q)$ is their Riemannian distance. Under these conditions it is easy to check that
$$k_1^*(t)={\omega_{n-1}\over n}(c_ d)^{-p'}t^{-1/p'}\Big(1+O\big(t^{-1/p'+\e}\big)\Big)$$
for small $t$, and it is clear that the estimate $k_2^*(t)\le C t^{-n/(\lambda p')}$ would follow if the underlying  Borel measure $\nu$ satisfies $\nu\big(B(P,r)\big)\le Cr^\lambda$, for small geodesic balls $B(P,r)$. These facts, and similar ones for vector-valued operators, imply  inequalities  such as  those of Theorems 6 and 7, and the sharpness statements are proven in essentially the same manner.

It would also be possible to extend this theorem to general Lorentz-Sobolev space, in the same spirit as in [Ci1], with  suitable and  slightly more general versions of our Theorems 1,1' and 4, which for simplicity we only treated in the  $L^p$ setting.

Finally, we wish to remark that our proof of Theorem 6 is of a somewhat different nature than the  one given by Cianchi in [Ci1]. In the special case $ d=1$ Cianchi started by applying the Sobolev inequality
$$\|\Psi\|_{L^{\lambda p\over n-p}(\Omega,d\nu)}^{}\le C\|\nabla \Psi\|_{L^p(\Omega)}^{}\eqno(66)$$
for some suitable $p<n$, with $\Psi(x)=\exp\Big[\ds{n-p\over p}\,\omega_{n-1}^{1\over n-1}\, |u|^{n'}(x)\Big]-1$. Note that (66) holds on $W_0^{1,p}$ and it is a special case of a general Sobolev inequality derived by Adams [Ad2], [Ad3]. Next, the fact that $\nabla e^u=e^u \nabla u$ combined with clever use of rearrangement inequalities and O'Neil's/Adams-Garsia's lemmas allowed Cianchi to obtain the result. However, it is unclear to us how to efficiently apply this strategy to arbitrary order derivatives.

\vskip1em\eject
\noin{\bf 5.  Sharp Moser-Trudinger inequalities for general elliptic differential operators}
\bigskip

In this section we extend Adams' original inequality to general elliptic differential operators on $\Rn$. We then specialize to particular cases. Some of these results can perhaps be further extended to suitable elliptic pseudodifferential operators, operators on manifolds, or even non-elliptic operators, but we will not treat such cases here.  

The structure of the fundamental solution of general elliptic operators is best explained by use of the pseudodifferential calculus, which we now very briefly recall.
\smallskip\noin{\bf Note.} 
Throughout this section all  functions will be {\it complex-valued}, unless otherwise specified.
\smallskip
\def\psido{{$\Psi$DO} }
A pseudodifferential (\psido) operator of order $d\in\R$
on an open set $U\subset\Rn$ 
is an operator $P:C_c^\infty(U)\rightarrow C^\infty(\Rn)$ of type 
$$Pf(x)=\op(p)f(x):=\irn\!\irn e^{2\pi i(x-y)\cdot\xi}p(x,\xi)f(y)dyd\xi$$
where $p\in C^\infty(U\times\Rn)$, the symbol of $P$, satisfies
 $$|\p_x^\beta \p_\xi^\alpha p(x,\xi)|\le C_{ d,\beta,K} (1+|\xi|)^{d-|\alpha|},\qquad x\in K,\, \xi\in\Rn$$
 for any $K\subset U$ compact, and any multiindices $\alpha,\beta$ (where $C$ is independent of $x,\,\xi$), and where $\p^{\alpha}_x=\partial_{x_1}^{i_1},...\partial_{x_n}^{i_n}$, if $\alpha=(i_1,...,i_n)$.

A classical, or polyhomogeneous, $\Psi$DO of order $d$ is given by a symbol $p(x,\xi)$ such that

\smallskip
\noin i) there is a sequence of functions $p_{d-j}(x,\xi)\in C^\infty(U,\Rn)$ which are homogeneous of order $d-j$ in $\xi$, for $|\xi|\ge1$:
 $$p_{d-j}(x,t\xi)=t^{d-j}p_{d-
j}(x,\xi), \qquad t\ge1,\;\;|\xi|\ge1,\;\;j=0,1,2... $$

\smallskip\noin ii) $p$ has the asymptotic expansion
$p\sim \sum_{j=0}^\infty p_{d-j},$ that is
$$\Big|\p_x^\beta \p_\xi^\alpha \Big[p(x,\xi)-\sum_{j=0}^{N-1} p_{d-
j}(x,\xi)\Big]\Big|\le C_{\alpha,\beta,K,N} (1+|\xi|)^{d-|\alpha|-N},\qquad x\in K,\,\xi\in\Rn$$
for any $M\in \N$, any multiindices $\alpha,\beta$, and any compact set $K\subset U$.

The principal symbol of such $P$ is the function $p_d(x,\xi)$ and the strictly homogeneous symbol
is defined by 
$$p_d^0(x,\xi)=|\xi|^d p_d(x,{\xi/|\xi|}),\qquad \xi\neq0,$$
that is the unique function on $\R^n\setminus \{0\}$ which coincides with $p_d$ for $|\xi|\ge1$, and which is homogeneous
or order $d$ in $\xi\in \Rn\setminus\{0\}$.

From now on we will consider only classical  \psido$\!\!$.
Every \psido can be written as an integral operator with kernel
$$K_P(x,y)=\irn e^{2\pi i(x-y)\cdot\xi}p(x,\xi)d\xi$$
defined in the sense of oscillatory integrals, and $C^\infty$ off the diagonal. If $K_P$ is $C^\infty(U\times U)$ then $K_P$ is called  smoothing (or negligible) operator.

For a classical \psido $T$ of order $d$, with $-n\le d< 0$  the  Schwarz kernel of $T$ has an expansion 
$$K(x,y)=\mathop{\sum}\limits_{0\le j<n-d}  k_{d-n+j}(x,x-y)+c(x)\log|x-y|+O(1)\eqno(67)$$
with $k_{d-n-j}(x,z)$ homogeneous of order $d-n-j$ in $z\in\Rn\setminus\{0\}$, $c(x)$ continuous on $U$, and $O(1)$ continuous and bounded on $K\times K$ for any $K\subset U$ compact. This fact is standard (see e.g. [Ca], Theorem 28) and follows  from taking the inverse   Fourier transform of the expansion of the symbol of $T$.
To be more specific,  let  $T$ have symbol $q$ with
$$q\sim \sum_{j=0}^\infty q_{-d-j}\quad x\in U,\, |\xi|\ge1,\eqno(68)$$
where $q_{-d-j}$ are homogeneous of order $-d-j$ for $|\xi|\ge 1$, and let  
$\F u=\hat u$ denote  the Fourier transform  on tempered distributions. Then for  $j<n-d$ the $q_{-d-j}^0$ are integrable around the origin, so that  modulo a smooth function ${\cal F}^{-1}\big[q_{-d-j}^0(x,\cdot)\big](z)$ coincides with
$$k_{d-n-j}(x,z):={\cal F}^{-1} \big[q_{-d-j}^0(x,\cdot)\big](z)\eqno(69)$$
which  is a homogeneous function of order $d-n-j$ in $z$, smooth in $x$ and in $z\neq 0$. When $j=n-d$ (which can only happen for integer $d$), then one can show (see [Ca], Theorem 27) that 
$${\cal F}^{-1} \big[q_{-n}(x,\xi)](z)=c(x)\log{1\over|z|}+O(1)\eqno(70)$$
as $z\to0$, with $c$ continuous, and in fact $c(x)=\ds{\int_{S^{n-1}}} q_{-n}^0(x,\omega)d\omega$ (the term $O(1)$ is in fact the sum of a  $C^\infty(\Rn\setminus0)$ homogeneous function of degree $0$ and a polynomial). Finally, the inverse Fourier transform of the error term $q-\sum_{j\le n-d}q_{-d-j}$ is easily estimated to be $O(1)$.

A \psido  operator $P$ is elliptic, if
$$p_d(x,\xi)\neq0,\qquad x\in U,\,|\xi|\ge 1$$
i.e. $p_d^0(x,\xi)\neq 0$ for $\xi\neq0$.  For elliptic operators one can construct a parametrix $G$, i.e. a \psido of order $-d$ such that $GP-I=PG-I$ is smoothing. The  symbol $q(x,\xi)$ of $G$ has an expansion like (68), where the $q_{-d-j}$  are determined from the symbol of $P$. In particular
$$ q_{-d}(x,\xi)={1\over p_d(x,\xi)}, \quad x\in U,\, |\xi|\ge1$$
so that the  parametrix admits a kernel expansion (67), with leading term given by 
$$k_{d-n}(x,x-y)={\cal F}^{-1} \big[p_{-d}(x,\cdot)^{-1}\big](x-y)=|x-y|^{d-n} g\big(x,{y-x\over|y-x|}\big)\eqno(71)$$
with $$g(x,\omega)=\bigg({1\over p_d^0(x,\cdot)}\bigg)^\wedge (\omega),\qquad\omega\in S^{n-1}\eqno(72)$$
the restriction of $k_{d-n}(x,\cdot)$ (which is $C^\infty(\R-\setminus\{0\})$) to the sphere.

\medskip\smallskip\bigskip
\noin{\it  Sharp inequalities  for the potential case}\medskip
\medskip
The precise asymptotic expansion of the parametrix operator $G$ suggests the following general theorem:
\bigskip
\proclaim Theorem 8. Let $\Omega$ be open and bounded on $\Rn$,  and let $K:\Omega\times \Omega\to\bar \R$ be  measurable and such that
$$K(x,y)=g\Big(x,{y-x\over|y-x|}\Big)\,|x-y|^{ d-n}+O\big(|x-y|^{ d-n+\e}\big)\eqno(73)$$ 
where $g:\Omega\times S^{n-1}\to \R$ is measurable and bounded, and $\e>0$.  For $0< d<n$, $\,p=\ds{n\over d}$, $\,\ds{{1\over p}+{1\over p'}}=1$, let, for $x\in \Omega$ 
$$Tf(x)=\int_\Omega K(x,y)f(y)dy,\quad f\in L^p(\Omega).$$ Then there exists $C>0$ such that 
$$\int_{\Omega} \exp\bigg[A^{-1}\bigg({|Tf(x)|\over \|f\|_p}\bigg)^{p'}\bigg]\,dx\le C\eqno(74)$$
for all $f\in L^p(\Omega)$, with 
$$A={1\over n}\,\sup_{x\in \Omega}\int_{S^{n-1}}|g(x,\omega)|^{p'}d\omega.\eqno(75)$$ 
\eject
If the supremum in (75) is attained at some $x_0\in\Omega$,  if $g(\cdot,\omega)$ is H\"older continuous of order $\sigma\in (0,1]$ at $x_0$ uniformly w.r. to $\omega$, i.e. if 
$$|g(x,\omega)-g(x_0,\omega)|\le C|x-x_0|^\sigma \qquad |x-x_0|\le \delta,\;\omega\in S^{n-1},$$ and  if $g(x_0,\cdot)$ is  H\"older continuous of order $\sigma$ on $S^{n-1}$,
then the constant $A^{-1}$ in (74) is sharp. In particular, there is  a suitable sequence  $r_m\to 0$ such that if   $\,E_m=B(x_0, r_m)\subseteq \Omega$ and  $\Phi_m(y)=K(x_0,y)|K(x_0,y)|^{p'-2}\chi_{\Omega\setminus E_m}^{}(y)$, then $\Phi_m\in
 L^{p}$ and
$$\lim_{m\to\infty}\int_\Omega \exp \bigg[\alpha\,\bigg({|T\Phi_m|\over\|\Phi_m\|_{p}^{}}\bigg)^{p'}\,\bigg]dx=+\infty,
\qquad \forall \alpha>{1\over A}.$$

\par
\medskip
\noin {\bf Remarks.}\smallskip\noin {\bf 1.} Cohn and Lu were the first to consider Adams inequalities for potentials of simpler type $g(y/|y|)|y|^{d-n}$, and the analogous version on the Heisenberg group ([CoLu1]).\smallskip \noin {\bf 2.} The H\"older continuity  condition on $g$ can be relaxed to an integral condition similar to that used in [CoLu1].
\medskip

 In view of Theorems 1 and 4 it is clear that to prove Theorem 8 it would essentially suffice to estimate the distribution function of the kernel $K$. This is done in the following lemma:
\proclaim Lemma 9. Suppose that $K$  is as in Theorem 8, satisfying (73) with $g$ bounded and measurable. Then for $s>0$ 
$$\sup_{x\in\Omega} \,|\{y\in\Omega: \,|K(x,y)|>s\}|\le   A s^{-p'}+O(s^{-p'-\sigma})\eqno(76)$$
 for suitable $\sigma>0$, with $A$ as in (75), with equality if  the sup in (75) is attained in $\Omega$. Moreover,   
$$\sup_{y\in\Omega}\,|\{x\in\Omega: \,|K(x,y)|>s\}|\le C  s^{-p'}.\eqno(77)$$
\par\smallskip
\noin {\bf Note.} A similar lemma was proved in [BFM], Lemma 2.3, for kernels in the CR sphere. 
\smallskip

\bigskip
\pf Proof of Lemma 9. From now on we will use the notation 
$$y^*={y\over|y|}.$$
 The hypothesis implies 
$$|K(x,y)|\le|g(x,(y-x)^*)|\,|x-y|^{ d-n}+C|x-y|^{ d-n+\e}$$
so that for any $x\in \Omega$
$$m_x(s):=|\{y\in\Omega: \,|K(x,y)|>s\}|\le |\{y\in\Rn: \,|g(x,y^*)|\,|y|^{ d-n}+C|y|^{ d-n+\e}>s\}|$$
and since
$$|g(x,y^*)|\,|y|^{ d-n}+C|y|^{ d-n+\e}>s\;\;\Longrightarrow \;\; |y|\le s^{-p'/n}\big(|g(x,y^*)|+C|y|^\e\big)^{p'/n}\le C s^{-p'/n}$$
then
$$m_x(s)\le{s^{-p'}\over n}\int_{S^{n-1}}\big(|g(x,y^*)|+Cs^{-\e p'/n}\big)^{p'}dy^*$$
which implies (76). 
Suppose that for some $x_0\in \Omega$
$$A={1\over n}\int_{S^{n-1}} |g(x_0,\omega)|^{p'}d\omega$$
and WLOG we can assume that $x_0=0$.
Since 
$$|K(0,y)|\ge|g(0,y^*)|\,|y|^{ d-n}-D|y|^{ d-n+\e}$$
for some $D>0$, then
$$\eqalign{m_0(s)&\ge |\{y\in\Omega: \,|g(0,y^*)|\,|y|^{ d-n}-D|y|^{ d-n+\e}>s\}|\cr&=|\{y\in\Omega: \,|y|<s^{-p'/n}\big(|g(0,y^*)|-D|y|^{\e}\big)^{p'/n}\}|.\cr}$$
But
$$|y|<s^{-p'/n}\big(|g(0,y^*)|-D|g(0,y^*)|^{\e p'/n}s^{-\e p'/n}\big)^{p'/n}\Longrightarrow |y|<s^{-p'/n}\big(|g(0,y^*)|-D|y|^{\e}\big)^{p'/n}\le C s^{-p'/n}$$
and if  $B(0,\delta)\subseteq \Omega$ then pick $s$ so large that $C s^{-p'/n}<\delta$. Let $\e<n/p'$ and  
$$E_s=\{y^*\in S^{n-1}:\,|g(0,y^*)|>D|g(0,y^*)|^{\e p'/n}s^{-\e p'/n}\}=\{y^*\in S^{n-1}:\,|g(0,y^*)|>D^{n\over n-\e p'}s^{-{\e p'\over n-\e p'}}\}.$$
Then,
$$\eqalign{m_0(s)&\ge {s^{-p'}\over n}\int_{E_s}\big(|g(0,y^*)|-D|g(0,y^*)|^{\e p'/n}s^{-\e p'/n}\big)^{p'}dy^*\ge {s^{-p'}\over n}\int_{E_s}\big(|g(0,y^*)|^{p'} -Cs^{-\sigma}\big)dy^*\cr&
\ge {s^{-p'}\over n}\int_{S^{n-1}}|g(0,y^*)|^{p'}dy^* -Cs^{-p'-\sigma}\cr} $$
which means that we have equality in (76).
Finally, (77) is a simple consequence of (73) and the boundedness of $g$.\endpf
\vskip-1em\eject
\pf Proof of Theorem 8. The previous lemma implies that 
$$K_1^*(t)\le At^{-1/p'}\big(1+O(t^{\e})\big),\qquad K_2^*(t)\le Ct^{-1/p'}\eqno(78)$$
so that the exponential inequality (74) follows form Theorem 1.

To prove sharpness, we appeal to Theorem 4. If the sup in  (75) is attained in $\Omega$, say WLOG at $x=0$,  then we have equality in the first estimate of (78). Choose $x_m=0$, and let $C_0,\, m,\, R$ large enough so that 
$$ \{y\in\Omega: \,|K(0,y)|>m\}\subseteq B(0,C_0m^{-p'/n})\subseteq \Omega\subseteq B(0,R).$$
Choosing 
$$r_m=C_0m^{-p'/n},\;\; E_m=B(0,r_m),\;\; B_m=B\big(0,\ts{1\over2} r_m\big)$$
we have that conditions a), b), c) of Theorem 4 are satisfies, so all we need to check is 
$$\int_{\Omega\setminus E_m} |K(x,y)-K(0,y)|\, |K(0,y)|^{p'-1} dy\le C\,,\qquad \forall x\in B_m.\eqno(79)$$
It is enough to verify this for $K(x,y)=g(x,(y-x)^*)|x-y|^{ d-n}$. By adding and subtracting $g\big(x,(y-x)^*\big)|y|^{ d-n}$ we see that  it suffices to verify
$$\int_{r_m\le |y|\le R} \big||x-y|^{ d-n}-|y|^{ d-n}\big|\,|y|^{- d}dy\le C,\qquad |x|\le {r_m\over2}\eqno(80)$$
which is the same as (62), and
$$\int_{r_m\le |y|\le R} |g\big(x,(y-x)^*\big)-g(0,y^*)|\,|y|^{-n}dy\le C,\qquad |x|\le {r_m\over2}.\eqno(81)$$
The H\"older continuity hypothesis on $g$ imply 
$$|g\big(x,(y-x)^*\big)-g(0,y^*)|\le C|x|^{\sigma}+C\bigg|{y-x\over|y-x|}-{y\over|y|}\bigg|^{\sigma}\le C|x|^{\sigma/2}|x-y|^{-\sigma/2}$$
but if $|x|\le r_m/2 $ and $|y|\ge r_m$, then $|x-y|\ge |y|/2$
and we are reduced to 
$$\int_{r_m\le |y|\le R} |x|^{\sigma/2}|y|^{-n-\sigma/2}dy\le C,\qquad |x|\le {r_m\over2}$$
which is clearly true.\endpf

\bigskip\eject
\noin{\it Sharp inequalities for  general elliptic operators}\medskip
With Theorem 8 at our disposal we are now in a position to extend Adams inequality (2) to rather general elliptic differential operators  of order 
$d<n$. In order for this machinery  to work, we need to make sure that we can write $u=T(Pu)$, for $u$ in a suitable class of smooth functions with compact support where $Pu\neq0$, and where $T$ is an integral operator with a kernel essentially equal to the  kernel of the parametrix $G$ of $P$. Note that we cannot simply take $T=G$, since $G$ is not, in general, an exact left inverse of P.

So let $P$ be a general  elliptic differential operator of  order $d$, written as 
$$P=\sum_{|\alpha|\le d} a_\alpha\p^\alpha$$
on the space of distributions $\D'(U)$, some open set $U$. The coefficients $a_\alpha$ are assumed to be $C^\infty$, complex-valued, and 
the principal symbol of $P$ then satisfies $p_d(x,\xi)=p_d^0(x,\xi)=(2\pi i)^d\sum_{|\alpha|=d} a_\alpha(x)\xi^\alpha\neq0$, for $\xi\neq0$. The adjoint of $P$ is the operator
$P^*=\sum_{|\alpha|\le d} (-1)^{|\alpha|}\p^\alpha\overline{a}_\alpha.$
\smallskip\medskip
 
\proclaim Theorem 10. Let $P$ be an elliptic,  differential operator of order $d<n$ on an open set $U$, with principal symbol $p_d(x,\xi)$. Let $\Omega$ be open and bounded with $\bar\Omega\subseteq U$, and let $p=\ds{n\over d}$, $\ds{1\over p}+\ds{1\over p'}=1$. If $P$ is injective  on $C_c^\infty(\bar\Omega)$, then there exists a constant $C$ such that 
$$\int_{\Omega} \exp\bigg[A^{-1}\bigg({|u(x)|\over \|Pu\|_p}\bigg)^{p'}\bigg]\,dx\le C\eqno(82)$$
for all $u\in W_0^{d,p}(\Omega)$, with 
$$A={1\over n}\,\sup_{x\in \Omega}\int_{S^{n-1}}|g(x,\omega)|^{p'}d\omega\eqno(83)$$ 
$$g(x,\omega)=\bigg({1\over p_d^0(x,\cdot)}\bigg)^\wedge (\omega),\qquad\omega\in S^{n-1}.\eqno(84)$$
In the special case  $p=p'=2$, i.e. $d=\ds{n\over2}$,  we have 
$$A={1\over n}\sup_{x\in\Omega} \int_{S^{n-1}}{1\over \big|p_{n/2}^0(x,\omega)\big|^2}d\omega=\sup_{x\in \Omega}\irn e^{-|p_{n/2}^0(x,\xi)|^2}d\xi.\eqno(85)$$
 If the supremum in (83) is attained in $\Omega$, and if the adjoint $P^*$ is injective on $C_c^\infty(\bar\Omega)$ (in particular if $P$ is self-adjoint) then the constant $A^{-1}$ in (82) is sharp.
\par 
\def\F{{\cal F}}
\pf Proof. It is enough to assume $u\in C_c^\infty(\Omega)$.  The given hypothesis assure the existence a (properly supported) \psido $T$ of order $-d$ such that $TP u=u$ for any $u\in C_c^\infty(\Omega)$ (in fact for any distribution $u$  with support in $\bar\Omega$) (see for example [Ca], Thm 24 and Thm 29). If $G$ denotes the parametrix of $P$ (in $U$) then $T=G+R$, some smoothing operator $R$, and therefore
$$u(x)=\int_\Omega K(x,y)P u(y)dy$$
with $K(x,y)$ having an expansion as in (67); in particular
$$K(x,y)=g\Big(x,{y-x\over|y-x|}\Big)\,|x-y|^{ d-n}+O\big(|x-y|^{ d-n+\e}\big)\qquad x,y\in \Omega$$ 
with 
$g(x,\omega)$ as in (84). Clearly this is precisely what one needs in order to apply Theorem 8, and hence prove   
 inequality (82). For the sharpness statement, if $P^*$ is injective on $C_c^\infty(\bar\Omega)$ then $P(Tf)=f$ in $\Omega$, for any $f\in \D'(U)$ (see [Ca], Thm 24),  $P:W_0^{d,p}(\Omega)\to L^p(\Omega)$ has  left inverse $T$, and $T:L^p\to W^{d,p}(\Omega)$. Suppose WLOG that the supremum in (83) is attained at $x_0=0\in\Omega$ and $B(0,1)\subset\subset\Omega$. If $\Phi_m$ is the sequence in the sharpness statement of Theorem 8, let $u_m=\psi T\Phi_m$, 
with $\psi\in C_c^\infty(\Omega)$ and $\psi=1$ in $B(0,1)$. Then $u_m\in W^{d,p}_0(\Omega)$ and it's easy to check using the Leibniz rule that $P u_m=\psi \Phi_m+ S \Phi_m$, where $S$ is a pseudodifferential operator of order  at most $-1$. Hence, since $|\Phi_m(y)|\le C |y|^{-d}$ on $\Omega$, we have $|S\Phi_m(y)|\le C|y|^{-d+1}$ and a straightforward estimate shows
$$\|Pu_m\|_p=\|\Phi_m\|_p+O(1),$$
 so that  
$$\int_{\Omega} \exp\bigg[\alpha\bigg({|u_m(x)|\over \|Pu_m\|_p}\bigg)^{p'}\bigg]\,dx\ge \int_{B(0,1)} \exp\bigg[\alpha\bigg({|T \Phi_m(x)|\over \|\Phi_m\|_p+C}\bigg)^{p'}\bigg]\,dx\to\infty$$
for any $\alpha>A^{-1}$.
Finally, when $p=2$ the first formula for $A$ given in (85) is a consequence of  the following  spherical Parseval's formula:     if $f,g\in C^\infty(S^{n-1})$ and $E_{-d} (f),\, E_{d-n} (g)$ are their homogeneous  extensions to $\R^n\setminus 0$ of order $-d$ and $d-n$ respectively $(0<d<n)$, then 
$$\int_{S^{n-1}}E_{-d}^\wedge(f)\overline{E_{d-n}^\wedge(g)}=\int_{S^{n-1}} f\overline g.\eqno(86)$$
The formula was originally found by Koldobsky  in case $f$ and $g$ are powers of Minkowsky functionals of smooth star bodies (see [K] for a proof). In [Mil] the above version is proven for $f,g$ real-valued and even, but a small modification of the proof yields the more general result. The second identity in (85) is obtained by a polar coordinate change. \endpf
\def\U{{\cal U}}
\vskip-2.5em
\noin{\bf Remarks.}\smallskip
 \noin{\bf 1.} It is possible to extend slightly Theorem 10 to the case when $P$ does not have a trivial nullspace. Indeed, in the setup of Theorem 10, if $\ker P$ denotes the nullspace of $P$ among distributions which are supported in $\bar\Omega$,  then $\ker(P)$ consists of $C^\infty$ functions, and it is finite-dimensional. 
\smallskip\noin{\bf 2.} The argument for the sharpness statement above can clearly be used to settle the sharpness statement of Adams  original result, or that of  our Theorem 6, in the case of $\Delta^{d/2}$ for $d$ even.\smallskip
\noin{\bf 3.} The purpose of the second identity in (85) is that in many situations the exponential integral can actually be evaluated explicitly. See a related calculation in Corollary 14 below.\smallskip
\smallskip \noin{\bf 4.} Both Theorems 8 and 10 can be formulated under a Borel measure satisfying  condition (52), in the same spirit as in Theorems 6 and 7 - the changes are minimal. We kept the usual Euclidean measures to better emphasize the relations between the sharp constants and the operators.

\smallskip

We now specialize Theorem 10 to the second order case. Let us start with an elliptic operator 
$$P=\sum_{j,k=1}^n a_{jk} \partial_{jk}^2+\sum_{j=1}^n b_j \partial_j+ c\eqno(87)$$
with $a_{jk},b_j,c\in C^\infty(U)$, {\it real-valued}, $a_{jk}=a_{kj}$, and let $\A_x=\big(a_{jk}(x)\big)$, an $n\times n$ symmetric matrix, which we assume to be positive definite. As before we assume that $\Omega$ is bounded, open and
 $\bar\Omega\subseteq U$. \medskip

\proclaim Corollary 11. Suppose that $P$ is an elliptic operator as in (87), and injective  on $C_c^\infty(\bar \Omega)$. Then for $n>2$ there exists $C>0$ such that 
$$\int_\Omega \exp\bigg[n(n-2)^{n\over n-2}\omega_{n-1}^{2\over n-2}\,\inf_{x\in \Omega}(\det\A_x)^{{1\over n-2}}\bigg({|u(x)|\over \|Pu\|_{n/2}}\bigg)^{n\over n-2}\bigg]dx\le C\eqno(88)$$
for all $u\in W_0^{2,n/2}(\Omega)$.
If $\,\ds{\inf_{x\in \Omega}}\det\A_x$ is attained in $\Omega$ then the exponential constant in  (88) is sharp.\par
\medskip
\noin{\bf Note.} If $P$ is strongly  elliptic in $U$, the classical theory (e.g. [GT], Thm 8.9) guarantees that $P$ is certainly injective on $C_c^\infty(U)$ if $c\le 0$. \medskip
\pf Proof. All we need to do is apply Theorem 10 to the operator $P$, with 
$$g(x,\omega)=-{1\over (2\pi)^2} \irn {e^{-2\pi i \omega\cdot\xi}\over \xi^T\A_x\xi} \,d\xi$$
where $\xi^T$ denotes the transpose of the vector $\xi$ seen as a column vector. If $\lambda_1(x),...,\lambda_n(x)$ denote the positive eigenvalues of $\A_x$ and if ${\xi\over\sqrt\lambda}=\Big({\xi_1\over\sqrt\lambda_1},...,{\xi_n\over\sqrt\lambda_n}\Big)$,  then for some orthogonal matrix $R$
$$g(x,\omega)=-{1\over\sqrt {\det \A_x}}\irn {e^{-2\pi i R\omega\cdot{\xi\over\sqrt\lambda}}\over(2\pi)^2|\xi|^2} \,d\xi=-{c_2\over \sqrt {\det \A_x}}\,\bigg|{R\omega\over\sqrt \lambda}\bigg|^{2-n}$$
where $c_2=\ds{1\over (n-2)\omega_{n-1}}$ is the constant in the Newtonian potential, as in (54).
\smallskip
Next, we compute 
$$\int_{S^{n-1}} |g(x,\omega)|^{n\over n-2}d\omega=\bigg({c_2\over\sqrt {\det \A_x}}\bigg)^{n\over n-2} \int_{S^{n-1}}\bigg|{\omega\over\sqrt \lambda}\bigg|^{-n}d\omega.$$
But the computations of the volume (with $x=x^*|x|$)
$$\eqalign{{\omega_{n-1}\over n}\,\sqrt {\det \A_x}&=\Big|\Big\{x:\,\Big|{x\over\sqrt\lambda}\Big|<1\Big\}\Big|=\Big|\Big\{x:\,|x|<\Big|{x^*\over\sqrt\lambda}\Big|^{-1}\Big\}\Big|={1\over n}\,\int_{S^{n-1}}\bigg|{\omega\over\sqrt \lambda}\bigg|^{-n}d\omega\cr}$$
give that 
$\int_{S^{n-1}}\big|{\omega/\sqrt \lambda}\big|^{-n}d\omega=\omega_{n-1}\sqrt {\det \A_x}$, and this concludes the proof.\hskip4em ///\smallskip\smallskip

\noin{\bf Remarks.}\smallskip\noin 
 \noin{\bf 1.} In case $b_1=...b_n=c=0$ the result of Corollary 11 can be derived directly from the known asymptotic expansion of the fundamental solution of $P$, and under even less restrictive smoothness conditions on the coefficients. In the case of $\lambda$-H\"older continuous coefficients  ($0<\lambda<1$) a classical result (see [Mi], Thm 19, VIII) guarantees that the equation $Pu=0$ has a fundamental solution $K(x,y)$ with an expansion
$$K(x,y)={c_2\over \sqrt {\det \A_x}}\Big((x-y)^T\A_x^{-1}(x-y)\Big)^{2-n\over2}\big(1+O(|x-y|^{\lambda})\big).$$
This expansion can also be extended to Dini-continuous coefficients or even under  weaker conditions [MMcO]. 
With the aid of such expansion the calculation of the distribution function of $K$ is straightforward, and produces the same constant as that of the above corollary. For the sharpness result, one just needs to make sure that estimate d) of Theorem 4 is verified, under milder smoothness conditions on the coefficients (and ultimately of the function $g(x,\omega)$).
 \smallskip
\noin{\bf 2.} In [FFV], Thm. 3.5,  an estimate such as (88) is derived  using a different method,  and for elliptic operators with much more general coefficients; the constant produced there is  $n(n-2)^{n\over n-2}\omega_{n-1}^{2\over n-2}$, under the ellipticity hypothesis $\xi^T\A_x\xi\ge |\xi|^2$. In such hypothesis and with smoother coefficients,  it is clear that our constant is in general greater (i.e. better), since $\det\A_x\ge1$. 

\medskip\noin {\it Sharp inequalities for vector-valued operators.}\medskip
We now offer a version  of Theorem 10 for vector-valued differential operators of type
$$\P=(P_j),\qquad P_j=\sum_{|\alpha|\le d} a_{j\alpha} \p^\alpha,\qquad j=1,2,...,\ell,\quad \ell\in\N\eqno(89)$$
with $a_{j\alpha}\in C^\infty$ and complex-valued, with  sharp statements  in the special case $p=2$, i.e.  $d=n/2$.
 
The goal is clearly to extend Adams' inequality for  the operators $\nabla\Delta^{d-1\over2}$ with $d$ odd, by mimicking the integration by parts that leads to the representation formula (57). For the scalar case one can represent $u$ in terms of $Pu$ essentially in a unique way, if $P$ is elliptic and injective;  in the vector-valued situation, on the other hand,   a question of ``optimal representation" of $u$ in terms  of $\P u$ arises, in order to obtain sharpness. The basic idea is to start with a vector-valued differential operator $\P$ as above,  and assume that for a given operator  ${\bf Q}=(Q_j)$ of order $d'$,  the operator $L={\bf Q}^*\cdot \P$ with order $d+d'\le n$ is elliptic and injective in $C^\infty_c$, so that it has an inverse $T$ of order $-d-d'$, and a Schwarz kernel $k(x,y)$. One can therefore write $u=(T{\bf Q}^*)\cdot \P u$ and apply Theorem 1' to obtain an Adams inequality, with exponential constant given explicitly in terms of the symbols of ${\bf Q}$ and $\P$. Clearly one cannot expect such constant to be sharp, given the dependence on $\bf Q$. We will not state in full generality such result, and for simplicity we will only deal with the case $\bf Q=\P$, since in the special situation  
$p=2$ i.e. $d=n/2$ a sharpness result can be easily obtained.

\def\Y{{\bf Y}}

For   vectors $\X=(X_j),\,\Y=(Y_j)$ we let $\X\cdot \Y=\sum_{j=1}^\ell X_j Y_j$, $\,|\X|=\big(\X\cdot\bar \X\big)^{1/2}=\Big(\sum_1^\ell|X_j|^2\Big)^{1/2}$.

\proclaim Theorem 12. Let $\P=(P_j)$ be an operator as in (89), with $d\le \ds {n\over2}$, defined on $\D'(U)$, some open set $U$. If $\Omega$ is open and bounded with $\bar\Omega\subseteq U$ and if  $L=\sum_1^\ell P_j^* P_j$
is elliptic on  $U$  and   injective  on $C_c^\infty(\bar\Omega)$, then there exists a constant $C$ such that, with $p=n/d$, 
$$\int_{\Omega} \exp\bigg[A^{-1}\bigg({|u(x)|\over \|\P u\|_p}\bigg)^{p'}\bigg]\,dx\le C\eqno(90)$$
for all $u\in W_0^{d,p}(\Omega)$, with
 $$A={1\over n}\,\sup_{x\in \Omega}\int_{S^{n-1}}|\g(x,\omega)|^{p'}d\omega$$ 
$$\g(x,z)=\big(g_j(x,z)\big),\qquad g_j^{}(x,z)=\Bigg({\bar p_j^0(x,\cdot)\over \ds{\sum_{k=1}^\ell |p_k^0(x,\cdot)|^2}}\Bigg)^\wedge(z),$$ 
where $p_j^0(x,\xi)=(2\pi i)^{d}\sum_{|\alpha|=d}a_{j\alpha}(x)\xi^\alpha$ is the principal symbol of $P_j$.\smallskip
In the case $p=2$ i.e. $d=\ds{n\over2}$ we have 
$$A={1\over n}\,\sup_{x\in \Omega}\int_{S^{n-1}}\bigg(\sum_{j=1}^\ell |p_j^0(x,\omega)|^2\bigg)^{-1}d\omega=\sup_{x\in \Omega}\irn \exp\bigg(-\!\!\sum_{j=1}^\ell |p_j^0(x,\xi)|^2\bigg)d\xi\eqno(91)$$ 
 and if the supremum in (91) is attained in $\Omega$, then the constant $A^{-1}$ in (90) is sharp.
\par
\medskip
\pf Proof. The given hypothesis on $L$ imply, just  as before, that  we can write any $u\in C_c^\infty(\Omega)$ as
$u=T(Lu)=\sum_{j} TP_j^*(P_j u)$, for a certain \psido $T$ of order $-2d\ge -n$, with Schwarz kernel $k(x,y)$ and principal symbol $p(x,\xi)=\Big(\sum_{k=1}^\ell |p_k^0(x,\xi)|^2\Big)^{-1}$. Since now $TP_j^*$ is a \psido
of order $-d$, and with principal symbol $\bar{p}_j^0(x,\xi)p(x,\xi)$, then it has a Schwarz kernel $K_j(x,y)$
so that $$K_j(x,y)=g_j^{}\big(x,(y-x)^*\big)|x-y|^{d-n}+O(|x-y|^{d-n+\e}).$$
The inequality in (90)  follows now from Theorem 1', since \def\K{{\bf K}} if $\K=(K_j)$ then 
$$u=\int_\Omega \K(x,y)\cdot \P u(y)dy$$
with $|\K(x,y)|=\big|\g\big(x,(y-x)^*\big)\big|\,|x-y|^{d-n}+O(|x-y|^{d-n+\e})$
and the estimates on its distribution functions follow from  Lemma 9. The formula for $A$ given in (91) is a consequence  of the spherical Parseval formula (86).

 To prove sharpness of the constant in (90) in the special case $p=2$,  we proceed as in the  proof of Theorem 6. Let the supremum in (91) be achieved at some $x_0\in \Omega$ and WLOG assume $x_0=0$. Note that $K_j(x,\cdot)=P_j k(x,\cdot)$, where $k$ is the kernel of $T$, and that $k(0,y)=c\log{1\over|y|}+O(1)$, some $c>0$, as per (70); let's say that 
$$c\log{c_0\over|y|}\le k(0,y)\le c\log{c_1\over|y|},\qquad y\in \bar\Omega$$
for some $c_0,c_1>0$. Now, using the same $\varphi$ as in (65), with $r_m\to0^+$ to be selected later, define
$$u_m(y)=\cases{0 &  for $\;k(0,y)\le \delta$\cr \cr\varphi\big(k(0,y)\big) & for $\;\delta< k(0,y)\le 1+\delta$\cr  k(0,y) & for $\;1+\delta< k(0,y)\le c \log\ds{1\over r_m}-1-\delta$\cr   c\log\ds{1\over r_m}-\varphi\Big({c\log\ds{1\over r_m}-k(0,y)\Big)} & for $c\log\ds{1\over r_m}-1-\delta<k(0,y)\le c\log\ds{1\over r_m}-\delta$\cr 
c\log\ds{1\over r_m}& for $\;k(0,y)>c\log\ds{1\over r_m}-\delta$.\cr}$$
 Then $u_m=0$ if $|y|>c_1 e^{-\delta/c}$, hence  we can choose $\delta$ so  large that the support of $u_m$ is inside $\Omega$, which implies that   $u_m\in W^{n/2,2}_0(\Omega)$. Additionally,  $u_m=c\log{1\over r_m}$ for $|y|<c_0 r_m e^{-\delta/c}$, and $u_m(y)=k(0,y)$ for $c_1r_m e^{1+\delta\over c}<|y|< c_0 e^{-{1+\delta\over c}}$, for $m$ large enough. So 
$$P_j u_m(y)=\cases{ K_j(0,y) & for $c_1r_m e^{1+\delta\over c}< |y|< c_0 e^{-{1+\delta\over c}}$\cr O(1) |y|^{-n/2}  
 & for $c_0 r_m e^{-\delta/c}< |y|<c_1r_m e^{1+\delta\over c}$\cr O(1) & otherwise\cr}$$
(here we used the chain and product rules, combined with $|\partial_y^\alpha k(0,y)|\le C |y|^{-|\alpha|}$, for $|\alpha|>0$, since $\partial^\alpha k$ is the kernel of the operator $\partial^\alpha T$, which has order $|\alpha|-n$).
\eject
Now  choose $r_m$ so that $\{y\in \Omega : |\K(0,y)|>m\}\subseteq B(0,C m^{-2/n})\subseteq B(0, c_1 r_m e^{1+\delta\over c})$, and therefore, we can apply (45) of Theorem 4 with $E_m=B(0, c_1 r_me^{1+\delta\over c})$ to conclude
$$\int_{\Omega\setminus B(0, c_1 r_me^{(1+\delta)/ c})} |\K(0,y)|^2dy=A\log{1\over r_m}+O(1)$$
which allows us to conclude $\|\P u_m\|_2^2=A\log{1\over r_m}+O(1)$ 
and the sharpness of the exponential constant follows immediately from (49), just as in the proof of Theorem 6. \endpf

We will give one first application of the above theorem to first order operators. Consider a family of  operators 
$$\P=\big(P_j\big)_{j=1}^n,\quad P_j=\sum_{k=1}^n a_{jk} \partial_k+b_j\eqno(92)$$
with  $a_{jk},\, b_j$  real-valued and $C^\infty$  on some open set $U\supseteq \bar\Omega$, with $\Omega$ bounded.

\proclaim Corollary 13. Suppose that $\A_x=\big(a_{jk}(x)\big)$ is invertible on $U$ and that  $L=\sum_{j=1}^n P_j^*P_j$ is injective  on $C_c^\infty(\bar \Omega)$. Then, for $n>1$ there exists $C>0$ such that 
$$\int_\Omega \exp\bigg[n\omega_{n-1}^{1\over n-1}\,\inf_{x\in \Omega}|\det\A_x|^{{1\over n-1}}\bigg({|u(x)|\over \|\P u\|_{n}}\bigg)^{n\over n-1}\bigg]dx\le C\eqno(93)$$
for all $u\in W_0^{1,n}(\Omega)$.
If  $\,\ds{\inf_{x\in \Omega}}|\det\A_x|$ is attained in $\Omega$ then the exponential constant in  (93) is sharp.\par
\bigskip
\noin{\bf Note.} In [FFV], Theorem 3.3, a similar estimate is given for less regular  coefficients,  under the condition $\xi^T\A_x\xi\ge |\xi|^2$, and with exponential  constant $n\omega_{n-1}^{1\over n-1}$, which is smaller than the one given in the above corollary.\smallskip

\medskip\pf Proof. The proof of (93) is just an application of Theorem 12. One just has to first compute $\g$, proceeding like in the proof of Theorem 12: if $\P_0=\Big(\sum_j a_{ij}\p_j\Big)=\A_x\cdot\nabla$
$$\eqalign{-\g(x,z)&=\bigg({(2\pi i)\A_x\xi\over(2\pi)^2 |\A_x\xi|^2}\bigg)^\wedge(z)= \P_0\bigg({1\over(2\pi)^2 |\A_x\xi|^2}\bigg)^\wedge(z)={1\over |\det \A_x|}\P_0\bigg({1\over (2\pi)^2|\xi|^{2}}\bigg)^\wedge\big((\A_x^{-1})^Tz\big)\cr&={c_2\over |\det \A_x|}\P_0 \big|(\A_x^{-1})^T z\big|^{2-n}={(2-n)c_2\over |\det \A_x|} \big((\A_x^{-1})^Tz\big)\big|(\A_x^{-1})^T z\big|^{-n}\cr}$$
since if $\A_x^{-1}=(a_{jk}')$ then  
$\sum_j a_{ij}\partial_j\big|(\A_x^{-1})^Tz\big|^{2-n}=(2-n)\big|(\A_x^{-1})^Tz\big|^{-n}\sum_{j,k} a_{ij}a_{jk}'\big((\A_x^{-1})^T z\big)_k. $ Estimate (93) follows since
$$A={1\over n}\sup_x\int_{S^{n-1}}|\g(x,\omega)|^{n\over n-1}d\omega={1\over n}\sup_x{1\over (\omega_{n-1}|\det \A_x|)^{n\over n-1}} \int_{S^{n-1}}\big|(\A_x^{-1})^T\omega\big|^{-n}d\omega$$
and $\ds\int_{S^{n-1}}\big|(\A_x^{-1})^T\omega\big|^{-n}d\omega=\omega_{n-1}|\det \A_x|$. For the sharpness statement, suppose  WLOG that $\,\ds{\inf_{x\in \Omega}}|\det\A_x|$ is attained at  $x_0=0\in\Omega$ and that the ellipsoid $\{y:|\A_0^{-1}y|<1\}\subseteq \Omega$. Take any $r_m\downarrow 0$, $r_m<1$, and let
$$u_m=\cases{\log|\A_0^{-1}y|^{-1} & if $r_m<|\A_0^{-1}y|<1$\cr
\log r_m^{-1} & if $|\A_0^{-1}y|\le r_m$\cr
0 & if $|\A_0^{-1}y|\ge1.$\cr}$$
 Then $u_m\in W^{1,n}(\Omega)$, $\P u_m(y)=-(\A_0^{-1}y)|\A_0^{-1}y|^{-2}+O\big(\log|\A_0^{-1}y|^{-1}\big)$ if $r_m<|\A_0^{-1}y|<1$, and it's easy to check that 
$\|\P u_m\|_n^n=\omega_{n-1} |\det\A_0|\log{1\over r_m}+O(1)$. The result follows from (49), with $B_m=\{y: |\A_0^{-1}y|<r_m\}$.
\endpf

As a second and final quick application of Theorem 12 we consider in $\R^4$ the  second order operators 
$$\P_1=\big(\partial_{11}^2,\,\partial_{22}^2,\,\partial_{33}^2,\,\partial_{44}^2\big)\,
\quad\P_2=\big(\partial_{11}^2+\partial_{22}^2,\,\partial_{33}^2+\partial_{44}^2\big)\,\quad 
\P_3=\big(\partial_{11}^2+\partial_{22}^2+\partial_{33}^2,\,\partial_{44}^2\big)$$ 
\proclaim Corollary 14. Let $\Omega\subseteq \R^4$ be open and bounded. Then there exists $C>0$ such that for $j=1,2,3$
$$\int_{\Omega} \exp\bigg[\,B_j\,\bigg({|u(x)|\over\|\P_j u\|_2}\bigg)^2\bigg] dx\le C\eqno(94)$$
with
$$B_1={\pi^4\over \Gamma\big({5\over4}\big)^4},\qquad B_2=64 \pi,\qquad B_3={16\pi^{5/2}\over\Gamma\big({3\over4}\big)}$$
for any $u\in W_0^{2,2}(\Omega)$, and the constants $B_j$  are sharp.\par

Note that the constant $32\pi^2$ in the sharp inequality 
$$\int_{\Omega} \exp\bigg[32\pi^2\bigg({|u(x)|\over\|\Delta u\|_2}\bigg)^2\bigg] dx\le C$$
is bigger than  all of the constants  in (94), in fact $32\pi^2>B_3>B_2>B_1$; this is   consistent with $\|\Delta u\|_2\ge \|\P_3 u\|_2\ge \|\P_2u\|_2\ge\|\P_1u\|_2$,  which is  easily seen via Fourier transform.
\medskip
\pf Proof. We can apply Theorem 12, since the operator $L=\P_1^*\cdot\P_1=\sum_1^4 {\partial^4\over\partial x_j^4}$ is elliptic and injective on $C_c^\infty(\bar\Omega)$, and the same is true for $\P_2^*\cdot \P_2$ and $\P_3^*\cdot \P_3$. The computation of the constants follows easily from (91) and the identity
$$\int_{\sR^m} \exp\bigg[-\bigg(\sum_{j=1}^m x_j^2\bigg)^{p/2}\bigg] dx={2\pi^{m/2}\Gamma\big(1+{m\over p}\big)\over m\Gamma\big({m\over2}\big)}$$
valid for $m\in \N$ and $p>0$. Note that $B_1^{-1}$ is in fact the volume of the convex body $\Big\{x\in \R^4:\,\sum_1^4 x_j^4<1\Big\}$  (see for example [K]).
\endpf\eject
\bigskip
\noin {\bf 6.  Sharp Adams inequalities for  sums of weighted potentials.}
\medskip
\smallskip
 As another   illustration of how   Theorems 1 and 4 can be used, we offer an  extension of Adams' inequality (3) in a different direction:

\bigskip\proclaim Theorem 15.  Let $\Omega,\,\Omega'$ be  bounded open sets of $\Rn$,  $a_1,...,a_N\in\Rn$,$\,a_j\neq a_k,\, j\neq k$. Let $\U$ be a bounded open set of $\R^n\times\R^n$, with  $\Omega'\times\Omega\subset\subset \U$, and  let  $g_j:\bar\U\to\R$,  be H\"older continuous of order $\sigma_j\in(0,1],\, j=1,2,...N$. For $0< d<n$, $\,p=\ds{n\over d}$, $\,\ds{{1\over p}+{1\over p'}}=1$, let, for $x\in \Omega'$ and $y\in\Omega$,
$$K(x,y)=\sum_{j=1}^N g_j(x,y)|x+a_j-y|^{ d-n},\qquad Tf(x)=\int_\Omega K(x,y)f(y)dy,\quad f\in L^p(\Omega),$$
If $\g=(g_1,...g_N)$ and 
$$M(\g):=\sup\bigg\{\sum_{j=1}^N |g_j(x,x+a_j)|^{p'},\;x\in\Omega',\,(x,x+a_j)\in\bar\U,\,j=1,...N\bigg\}>0,\eqno(95)$$
then there exists $C$ such that
for any $f\in L^p(\Omega)$
$$\int_{\Omega'} \exp\bigg[ {n\over\omega_{n-1} M(\g)}\bigg({|Tf|\over \|f\|_{p}}\bigg)^{p'}\,\bigg]dx\le C\eqno(96)$$
with   $C$ independent of $f$. 
If 
$$\Omega^*:=\Omega'\cap\bigcap_{j=1}^N(\Omega-a_j)\neq\emptyset\eqno(97)$$
and $M(\g)$ is attained on $\Omega^*$, then the constant $\ds{n\over\omega_{n-1} M(\g)}$ is sharp in (96), i.e. it cannot be replaced by a larger constant.
\par
\medskip

\pf Proof.  Fix $x\in \Omega'$. If $\delta>0$ is such that $\delta<{\rm dist} \big(\bar{\Omega'}\times\bar\Omega\,,\,\bar \U^c\big)$, and   $B(a_j,\delta)\cap B(a_k,\delta)=\emptyset$,for $j\neq k$, $j,k=1,...N$,  then for $s>s_1:=N\delta^{ d-n}\max_j \|g_j\|_\infty^{}$ we have

$$\big|\{y\in\Omega:\,|K(x,y)|>s\,\}\big|=\sum_{j=1}^N\big|\{y\in\Omega\cap B(x+a_j,\delta):\, 
|K(x,y)|>s\}\big|$$

With our choice of $\delta$ it's clear that if $(x,x+a_j)\notin \bar\U$ then $\Omega\cap B(x+a_j,\delta)=\emptyset$ so 
$$\big|\{y\in\Omega\cap B(x+a_j,\delta):\, |K(x,y)|>s\}\big|=0$$
for any $s>s_1$ (in fact for  any $s>0$).

Assume that  $(x,x+a_j)\in\bar\U$ and $y\in \Omega\cap B(x+a_j,\delta)$. Then 

$$|K(x,y)|\le |g_j(x,y)|\,|x+a_j-y|^{ d-n}+C\delta^{ d-n}\le |g_j(x,x+a_j)|\,|x+a_j-y|^{ d-n}+C|x+a_j-y|^{ d-n+\epsilon}$$
some $\epsilon>0$, $\epsilon<n- d$ and $C$ independent on $x,y$. As a consequence, if $|K(x,y)|>s$ then
$$|x+a_j-y|<s^{-1/(n- d)}\big(|g_j(x,x+a_j)|+C|x+a_j-y|^\epsilon\big)^{1/(n- d)}\le Cs^{-1/(n- d)}$$

and 

$$\eqalignno{\big|\{y\in\Omega\cap B(x+a_j,\delta):\, 
|K(x,y)|>s\}\big|&\le {\omega_{n-1}\over n}\, s^{-n/(n- d)}\big(|g_j(x,x+a_j)|+Cs^{-\epsilon/(n- d)}\big)^{n/(n- d)}\cr&\le{\omega_{n-1}\over n}\, s^{-p'}|g_j(x,x+a_j)|^{p'}+Cs^{-p'-\sigma}&(98)\cr} $$

some $\sigma>0$ (we used here, for example, that $|(a+b)^\nu-b^\nu|\le Ca^{\min\{1,\nu\}}$ if $\nu>0$ and $a,b\in[0,K]$, some fixed $K>0$, C independent of $a,b$).

Now we see that if $x\in \Omega'$  and $(x,x+a_j)\in\bar\U$ for all $j$, then
$$\big|\{y\in\Omega:\,|K(x,y)|>s\,\}\big|\le s^{-p'}{\omega_{n-1}\over n}\,\sum_{j=1}^N |g_j(x,x+a_j)|^{p'}+O\big(s^{-p'-\sigma}\big),\quad \forall s> s_1\eqno(99)$$
(with $|O\big(s^{-p'-\sigma}\big)|\le Cs^{-p'-\sigma}$, $C$ independent of $x,s$), from which it follows that 

$$\sup_{x\in\Omega'} m\big(K(x,\cdot),s\big)\le s^{-p'}M(\g)+O\big(s^{-p'-\sigma}\big).\eqno(100)
$$

On the other hand, the same argument used to derive (99) can be used to show

$$\big|\{x\in\Omega':\,|K(x,y)|>s\,\}\big|\le B s^{-p'}\eqno(101)$$
for all $s>s_1$, and $y\in\Omega$, for some $B>0$ independent of $y$.  

Estimate (96) now follows from Theorem 1, using (100), (101) together with Fact 3.

Now assume that $\Omega^*\neq\emptyset$ and that the sup in (95) is attained inside $\Omega^*$, say at $x^*$. WLOG we can assume that $x^*=0$ (indeed it's enough to perform a translation by $x^*$ in both the $x$ and the $y$ variables). If $y\in \Omega\cap B(x+a_j,\delta)$
$$|K(0,y)|\ge |g_j(0,a_j)|\,|a_j-y|^{ d-n}-C|a_j-y|^{ d-n+\e}$$
and $|g_j(0,a_j)|\,|a_j-y|^{ d-n}-C|a_j-y|^{ d-n+\e}>s$ if and only if
$$|a_j-y|<s^{-1/(n- d)}\big(|g_j(0,a_j)|-C|a_j-y|^\epsilon\big)^{1/(n- d)}.
$$

Then letting 

$$\phi(s):=s^{-1/(n- d)}\big(|g_j(0,a_j)|-Cs^{-\e/(n- d)}|g_j(0,a_j)|^{\e/(n- d)}\big)^{1/(n- d)}\eqno(102)$$ we have 
$$\eqalign{&\big\{y\in\Omega\cap B(a_j,\delta):\,|a_j-y|<\phi(s)\big\}\cr&\subseteq \big\{y\in\Omega\cap B(a_j,\delta):\,|a_j-y|<s^{-1/(n- d)}\big(|g_j(0,a_j)|-C|a_j-y|^\epsilon\big)^{1/(n- d)}\big\}\cr
&\subseteq\big\{y\in\Omega\cap B(a_j,\delta):\,|K(0,y)|>s\big\}.\cr}\eqno(103)$$

Since $a_j\in\Omega$ let $\delta_0>0$ be such that $B(a_j,\delta_0)\subseteq\Omega$. There exists $s_0>s_1$ such that $0\le \phi(s)<\delta_0$ for all $s\ge s_0$, so that 

$$\big|\big\{y\in\Omega\cap B(a_j,\delta):\,|K(0,y)|>s\big\}\big|\ge {\omega_{n-1}\over n}\,\big(\phi(s)\big)^n\ge{\omega_{n-1}\over n}\, s^{-p'}|g_j(0,a_j)|^{p'}-Cs^{-p'-\sigma}$$
for all $s\ge s_0$.

This means that for all $s>s_0$
$$\big|\{y\in\Omega:\,|K(0,y)|>s\,\}\big|=s^{-p'}{\omega_{n-1}\over n}\,\sum_{j=1}^N |g_j(0,a_j)|^{p'}+O\big(s^{-p'-\sigma}\big)=s^{-p'}M(\g)+O\big(s^{-p'-\sigma}\big).$$

\smallskip\noindent

Now let us choose $x_m=0$ for $m\in\N$, 

$$E_m=\{y\in\Omega:\,|K(0,y)|>m\,\}=\bigcup_{j=1}^N\{y\in\Omega\cap B(a_j,\delta):\, 
|K(0,y)|>m\}$$
the union being disjoint for $m>s_1$. From  ii) we have $|E_m|=m^{-p'}M(\g)+O(m^{-p'
-\sigma})\to0$ as $m\to\infty$. Moreover, from (102) and (103), if $g(0,a_j)\neq0$ then $\{y\in\Omega\cap B(a_j,\delta):\, 
|K(0,y)|>m\}$ contains a ball of center $a_j$ and radius $C_j m^{-p'/n}$ some $C_j>0$, for all $m>m_j>s_1$; let $C_0$ be the smallest of such $C_j$ and let $$r_m=C_0 m^{-p'/n},\quad B_m=B\big(0,\ts{1\over2}r_m\big).\eqno(104)$$

With these choices conditions a), b), c)  of Theorem 4 are satisfied, so all we need is to check (43), i.e.
$$\int_{\Omega\setminus E_m} |K(x,y)-K(0,y)|\, |K(0,y)|^{p'-1} dy\le C\,,\qquad \forall x\in B_m\eqno(105)$$
some $C$ independent of $x$ and $m$.

Now observe the following elementary inequalities, valid for any $x\in \Omega'$ and $y\in \Omega$
$$\eqalign{|K(0,y)|^{p'-1}&\le C\sum_{j=1}^N |g_j(0,y)|^{ d/(n- d)}|y-a_j|^{- d}\cr&\le C\sum_{j=1}^N\big(|y-a_j|^{\e d/(n- d)}+|g(0,a_j)|^{ d/(n- d)}\big)|y-a_j|^{- d}\cr}\eqno(106)$$
$$\eqalign{\big|g_j(x,y)|x+a_j-y|^{ d-n}-g_j&(0,y)|a_j-y|^{ d-n}\big|\le C\big(|y-a_j|^\e+|x|^\e\big)|x+a_j-y|^{ d-n}\cr&\qquad +|g_j(0,a_j)|\big||x+a_j-y|^{ d-n}-|a_j-y|^{ d-n}\big|\cr}\eqno(107)
$$

$$\eqalign{|K(x,&y)-K(0,y)|\, |K(0,y)|^{p'-1}\le C\sum_{j=1}^N\bigg\{|x|^\e|y-a_j|^{- d}|x+a_j-y|^{ d-n}+\cr& +|g(0,a_j)|^{n/(n- d)}|y-a_j|^{- d}\big||x+a_j-y|^{ d-n}-|a_j-y|^{ d-n}\big|\bigg\}+\Phi(x,y)\cr}\eqno(108)$$
where $\Phi(x,y)\ge0 $ is integrable in $y\in B(0,R)$ some $R$ large enough so that $\int_{B(0,R)}\Phi(x,y)dy\le C$, independent of   $x\in \Omega'$.

By virtue of (108) it is enough to consider those $j$ for which $g(0,a_j)\neq0$, and for such $j$ we can write $\Omega\setminus E_m\subseteq \Omega\setminus B(a_j,r_m)$ (recall the definition of $r_m$ in (104)). Thus, it all boils down to (62), which we already checked, and the estimate
$$\sup_{|x|\le r_m/2}\;\int_{r_m\le|y|\le R} |x|^\e|y|^{- d}|x-y|^{ d-n}dy\le C,\eqno(109)$$
which is an easy consequence of (63).

\endpf

\bigskip

\noin{\bf Remarks.} \smallskip\item{\bf 1.} From the proof above it should be  apparent that Theorem 15 holds verbatim for kernels of type
$$K(x,y)=\sum_{j=1}^N g_j(x,y)|x+a_j-y|^{ d-n}\bigg[1+O\bigg(\sum_{j=1}^N |x+a_j-y|^{\e_j}\bigg)\bigg]\eqno(110)$$
where $\e_1,...\e_N>0$.
\smallskip
\item{\bf 2.} The regularity hypothesis on the $g_j$ can be somewhat relaxed to an integral
condition of type (43).\smallskip  
\item{\bf 3.} If the sup defining $M(\g)$ is not attained in $\Omega^*$, or if $\Omega^*$ is empty, then
the sharp constant in (96) will in general be larger, and the geometries of  the domains could play a definite role. For example, if $K(x,y)=|x-y|^{ d-n}$, and  $\Omega',\, \Omega$ are two open balls with empty intersection but tangent to one another (or two $C^1$ domains with the same property), then
$M(\g)=1$, but it's easy to see that the sharp constant in (96) is $2n/\omega_{n-1}$. This can bee seen by explicit asymptotics of the distribution function of the kernel with the given domains, together with Theorems 1 and 4. Similar considerations could be made if $\partial \Omega$ has corners, or even positive measure.
On the other hand, if $K(x,y)=|x+e_1-y|^{ d-n}+|x-e_1-y|^{ d-n}$, with $e_1=(1,0,..,0)$, and $\Omega'=B(0,10),\, \Omega=B\big(0,{1\over2}\big)$, then $\Omega^*=\emptyset$, $M(\g)=2$, but the sharp constant in (96) is $n/\omega_{n-1}$. This can be seen for example by splitting $\Omega'$ into two halves  each containing  $e_1$ or $-e_1$, and  noticing that in each half only one of the two potentials is really effective (i.e. Theorems 1 and 4 apply in each half separately).
\bigskip
On the $n-$dimensional Euclidean sphere $S^n$  Theorem 15 takes a somewhat  simpler form. Let $\eta,\xi$ denote points on $S^n$, and let $d\eta$ denote the standard volume element of $S^n$. 

\proclaim Theorem 16. On $S^n$ consider an operator 
$$Tf(\xi)=\isn K(\xi,\eta)f(\eta)d\eta,\qquad f\in L^1(S^n)$$
with $$K(\xi,\eta)=\sum_{j=1}^N g_j(\xi,\eta)|R_j\xi-\eta|^{ d-n}+O\Big(\sum_{j=1}^N |R_j\xi-\eta|^{ d-n+\epsilon_j}\Big),\quad 0< d<n,\;\sigma_j>0$$
for some $R_1,...,R_N\in SO(n)$, and $g_j:S^n\times S^n\to\R$ H\"older continuous of orders $\sigma_1,...,\sigma_N\in(0,1]$. If $p=\ds{n\over d},\,\ds{{1\over p}+{1\over p'}=1}$, $\,\g=(g_1,...,g_N^{})$ and if
$$M(\g)=\max_{\xi\in S^n}\sum_{j=1}^N |g_j(\xi,R_j\xi)|^{p'}>0,$$
 then there exists $C$ so that
$$\isn \exp\bigg[{n\over \omega_n M(\g)}\bigg({|Tf|\over\|f\|_p}\bigg)^{p'}\bigg]\,d\xi\le C\eqno(111)$$
for any $f\in L^p(S^n)$. The constant $\ds{n\over \omega_n M(\g)}$ in (111) is sharp.
\par
The proof of this theorem is identical to the one of Theorem 15, with the obvious modifications, and with the additional simplifications due to the compactness of $S^n$.\bigskip

On a compact Riemannian manifold $M$, with volume element $dV(P)$  and geodesic distance $d(P,Q)$, we have the following slight extension of Fontana's result ([F], Thm. 1.9):

\proclaim Theorem 17. On the compact Riemannian manifold $M$ consider an integral operator 
$$Tf(P)=\int_M K(P,Q)f(Q)dV(Q),\qquad f\in L^1(M)$$
with
$$K(P,Q)=g(P,Q)\,d(P,Q)^{ d-n}+O\big(d(P,Q)^{ d-n+\e}\big),\qquad 0< d<n,\;\e>0$$
with $g:M\times M\to\R$ H\"older continuous of order $\sigma\in(0,1].$ 
If $p=\ds{n\over d},\,\ds{{1\over p}+{1\over p'}=1}$ and if
$$M(g)=\max_{P\in M}|g(P,P)|^{p'}>0,$$
 then there exists $C$ so that
$$\int_M \exp\bigg[{n\over \omega_n M(g)}\bigg({|Tf|\over\|f\|_p}\bigg)^{p'}\bigg]\,dV(P)\le C\eqno(112)$$
for any $f\in L^p(M)$. The constant $\ds{n\over \omega_n M(g)}$ in (112) is sharp.
\par
The proof of Theorem 17 is a consequence of Theorems 1 and 4, and a sharp asymptotic estimate
of the distribution function of $K$, which is the same one as in the Euclidean case (Theorem 15)
given the fact  that the volume of a small geodesic ball is asymptotically the same as that of an Euclidean ball.
\bigskip
\noin{\bf 7. Sharp Adams inequalities on the CR sphere}\bigskip
As we mentioned in the introduction, Moser-Trudinger inequalities have recently been introduced in the context of CR-manifolds, first by Cohn and Lu [CoLu1,CoLu2] and more recently by Branson, Fontana, Morpurgo [BFM].  In [BFM], a special case of Theorem 1 of the present paper was quoted and used to derive sharp Adams inequalities for a class of convolution operators on the CR sphere ([BFM], Thm. 2.2.). The proof that such inequalities are sharp was only hinted in [BFM]; in this section we will provide a more detailed argument as an application of Theorem 4. 

We will now briefly recall 
the main setup.
Let $\Sn$ be the $(2n+1)-$dimensional sphere with its standard CR structure, i.e. that induced naturally from the ambient space $\C^{n+1}$, endowed with Hermitian product $\zbe$, where $\z=(\z_1,...,\zn)$, $\eta=(\eta_1,...,\eta_{n+1})$. The homogeneous dimension of $\Sn$ is denoted by $Q=2n+2$. Let $d\z$ be the standard volume element of the sphere, and $\omegan=2\pi^{n+1}/n!$ its volume; the average of a function $F$ on $\Sn$ is denoted by $\ds\avg F$. 

The Heisenberg group $\Hn$, with elements $(z,t)\in \C^n\times\R$ and group law
$(z,t)(z',t')=(z+z',t+t'+2\Im z\cdot \bar z')$
is biholomorphically equivalent to $\Sn$ via the Cayley transform $\Ca:\Hn\to S^{2n+1}\setminus(0,0,...,0,-1)$  given by 
$$\Ca(z,t)=\Big({2z\over 1+|z|^2+i t},{1-|z|^2-i t\over 1+|z|^2+i t}\Big)$$
and with inverse
$$\Ca^{-1}(\z)=\Big({\z_1\over1+\z_{n+1}},...,{\z_n\over1+\z_{n+1}},Im {1-\z_{n+1}\over
  1+\z_{n+1}}\Big).$$

The homogeneous norm on $\Hn$ is defined by 
$$|(z,t)|=(|z|^4+t^2)^{1/4}$$
and the distance from   $u=(z,t)$ and  $v=(z',t')$ is given as
$$d((z,t),(z',t')):=|v^{-1}u|=\big(|z-z'|^4+(t-t'-2\Im(z\bar
z'))^2\big)^{1/4}$$

On the sphere the distance function is defined as 
$$d(\z,\eta)^2:=2|1-\z\cdot \bar\eta|=\big|\,|\z-\eta|^2-2i
\,\Im(\z\cdot\bar\eta)\big|=
\big(|\z-\eta|^4+4\cdot{\rm \Im}^2(\z\cdot\bar\eta)\big)^{1/2}$$
and a simple calculation shows that if $u=(z,t),v=(z',t')$, and $\z=\Ca(u),\,\eta=\Ca(v)$.
then
$${|1-\z\cdot \bar\eta|\over
  2}=|v^{-1}u|^2\big((1+|z|^2)^2+t^2\big)^{-1/2}\big((1+|z'|^2)^2+(t')^2\big)^{-1/2}\eqno(113)$$
Further, we let 

$$u=(z,t)\in\Hn,\quad \Sigma=\{u\in\Hn:\,|u|=1\},\quad u^*={u\over|u|}=(z^*,t^*)\in\Sigma$$
$$\z=\Ca(u),\quad {1-\zn\over 1+\zn}=|z|^2+it=(|z|^4+t^2)^{1/2} e^{i\theta},\quad \n=\Ca(0,0)=(0,0,...,1),$$
and for $w\in\C,\,|w|<1$ we let
$$\theta(w)=\arg{1-w\over 1+w}\in \Big[-{\pi\over2},{\pi\over2}\Big].$$

 A function depending only on $\theta=\sin^{-1}t^*$ can be regarded as a function on the Heisenberg sphere $\Sigma$.
\medskip\eject
\proclaim {Theorem 18 ([BFM], Thm. 2.2)}. For  $0<d<Q$ let  $p=\ds{Q\over d}$ and  $\ds{1\over p}+{1\over p'}=1$.  Define
$$TF(\z)=\int_{\Sn} G(\z,\eta)F(\eta)d\eta,\qquad F\in L^p(\Sn)$$
where 
$$\eqalign{G(\z,\eta)&=g_0\big(\theta(\zbe)\big)\,d(\z,\eta)^{d-Q}+O\big(d(\z,\eta)^{d-Q+\e}\big)=\cr&=
2^{{d-Q\over2}} g_0\big(\theta(\zbe)\big)\,|1-\zbe|^{d-Q\over2}+O\big(|1-\zbe|^{{d-Q+\e\over2}}\big),\quad\z\neq\eta\cr}$$
for bounded and measurable $g_0:\big[\!-{\pi\over2},{\pi\over2}\big]\to\R$, with  
$\big|O\big(|1-\zbe|^{{d-Q+\e\over2}}\big)\big|\le C|1-\zbe|^{{d-Q+\e\over2}}$, some $\epsilon>0$, and with $C$ independent of $\z,\eta$.
\smallskip Then, there exists $C_0>0$ such that for all $F\in L^p(\Sn)$
$$\int_{\Sn}\exp\bigg[A_d\bigg({|TF|\over\|F\|_p}\bigg)^{p'}\bigg]d\z\le C_0\eqno(114)$$
with
$$A_d={2Q\over\ds{\int_\Sigma}|g_0|^{p'} du^*}.\eqno(115)$$
 Moreover, if the function $g_0(\theta)$  is H\"older continuous of order $\sigma\in(0,1]$ then the constant in (115) is sharp, in the sense that if it is replaced by a larger constant then there exists a sequence $F_m\in L^p(\Sn)$ such that the exponential integral in (114) diverges to $+\infty$ as $m\to\infty$. \par\smallskip

 In [CoLu1] Cohn and Lu give a similar result in the context of the Heisenberg group, and for kernels of type $G(u)=g(u^*)|u|^{d-Q}$, i.e. without any perturbations. An $\Hn$ version of Theorem 18 holds with virtually the same proof (in fact somewhat easier), but the two versions do not seem to be a consequence of each other.
\smallskip
In view of Theorem 1, to prove (114) it is enough to find an asymptotic estimate for  the distribution function of $G$. This is provided by the following result (which was proved in [BFM]):\smallskip 
\proclaim Proposition 19 {([BFM] Lemma 2.3)}. Let $G:\Sn\times\Sn\setminus\{(\zeta,\zeta),\z\in\Sn\}\to\R $,  be measurable and such that 
$$G(\z,\eta)=g\big(\theta(\zbe)\big)\,|1-\zbe|^{-\alpha}+O\big(|1-\zbe|^{-\alpha+\e}\big),\quad\z\neq\eta$$
 some bounded and measurable $g:\big[-{\pi\over2},{\pi\over2}\big]\to\R$, with  
$\big|O\big(|1-\zbe|^{-\alpha+\e}\big)\big|\le C|1-\zbe|^{-\alpha+\e}$, some $\epsilon>0$, and with $C$ independent of $\z,\eta$.
Then, for each $\eta\in \Sn$ and as $s\to+\infty$
$$\big|\{\z:\,|G(\zeta,\eta)|>s\}\big| = s^{-Q/2\alpha}\,{2^{Q/2-1}\over Q}\ints|g|^{Q/2\alpha}du^*+O\big(s^{-Q/2\alpha-\sigma}\big)$$
 for a suitable $\sigma>0$.
\par
Observe that $G$ above  may not be symmetric, but has upper and lower bounds with enough symmetries, so that in effect 
$G(\zeta,\cdot)^*(t)$ and $G(\cdot,\eta)^*(t)$, have the same asymptotic expansion in $t$ (independent of $\z,\eta$).
Proposition 19  combined with Theorem 1 gives (114). We now apply Theorem 4 in order to show the sharpness statement (this part was not done in [BFM]).
\smallskip
\pf Proof  of sharpness statement of Thm 18. The proof is similar to that of Theorem~8. Let $\z_m=\n$ , $r_m=C_0m^{-1/(Q-d)}<1$, so that 
$$\big\{\eta:\,|G(\n,\eta)|>m\big\}\subseteq E_m:=\big\{\eta:|1-\en|<2r_m^2\big\}$$
and let $B_m=\big\{\zeta:|1-\zn|<\ts{1\over4}r_m^2\big\}.$
\smallskip
 Conditions a), b), c) of Theorem 4 are met, from Proposition 19 and Remark 1 after Theorem 4, so all we need to do is show that 
$$\int_{\Sn\setminus E_m} |G(\zeta,\eta)-G(\n,\eta)|\, |G(\n,\eta)|^{p'-1}d\eta\le C\,,\qquad \forall \eta\in B_m.
$$

WLOG we can assume that $G(\zeta,\eta)=g\big(\theta(\zbe)\big)\,|1-\zbe|^{d-Q\over2}$, with $g=2^{d-Q\over2}g_0$; as it will be apparent from the proof below,  an
error term of type $|1-\zbe|^{{d-Q+\e\over2}}$ will produce an integrable function on $\Sn$, with uniformly bounded integral. So let us show that
$$\mathop\int\limits_{|1-\en|\ge 2r_m^2} |G(\zeta,\eta)-G(\n,\eta)|\, |G(\n,\eta)|^{d\over Q-d} d\eta\le C\,,\qquad  |1-\zn|< \ts{1\over4}r_m^2
\eqno(116)$$

By adding and subtracting the quantity \def\nbe{\n\cdot\bar\eta}
$g\big(\theta(\zbe)\big)|1-\en|^{d-Q\over2}$ we are reduced to proving the following estimates
$$\mathop\int\limits_{|1-\en|\ge 2r_m^2} \big|g\big(\theta(\zbe)\big)-g\big(\theta(\nbe)\big)\big||g(\nbe)|^{d\over Q-d}|1-\en|^{-Q/2}d\eta\le C\eqno(117)$$
$$\mathop\int\limits_{|1-\en|\ge2 r_m^2}\big|g\big(\theta(\zbe)\big)\big|\big|g\big(\theta(\nbe)\big)\big|^{d\over Q-d}\big||1-\zbe|^{d-Q\over2}-|1-\en|^{d-Q\over2}\big| \,|1-\en|^{-d/2}d\eta\le C\eqno(118)$$
valid for all $\zeta\in B_m$.

The first step is to transfer these integrals to $\Hn$ via the Cayley transform. Recall that the volume density of the Cayley transform is 
$$|J_\Ca(z,t)|={2^{2n+1}\over \big((1+|z|^2)^2+t^2\big)^{n+1}}\le {2^{2n+1}\over (1+|u|^4)^{Q/2}}$$

If $u=(z,t),\,v=(z',t')$, and $\z=\Ca(u),\,\eta=\Ca(v)$ then for $m$ large enough (using (113))
 $$2r_m^2\le|1-\en|={2|v|^2\over(1+2|z'|^2+|v|^4)^{1/2}}\le2|v|^2\;\;\Longrightarrow\;\; |v|\ge r_m$$
$${r_m^2\over4}>|1-\zn|= {2|u|^2\over(1+2|z|^2+|u|^4)^{1/2}}\ge {2|u|^2\over1+|u|^2}\;\;\Longrightarrow\;\; |u|<{r_m\over2}$$
and so if $\eta\notin E_m,\,\zeta\in B_m$ then (using that $|v^{-1}u|$ is a distance)
$${|1-\zbe|^{1/2}\over|1-\en|^{1/2}}={|v^{-1}u|\over|v|}\,{1\over(1+2|z|^2+|u|^4)^{1/4}}\ge{{1-\ds{|u|\over|v|}\over(1+|u|^2)^{1/2}}}\ge{1\over 2\sqrt 2}.$$
Since $|v^{-1}u|/|v|=1+O(|u|/|v|)$, for our range of $u$ and $v$, we obtain that the integrand in (118) in $\Hn$ coordinates is bounded above by
$$J_m=C\int_{|v|\ge r_m}\bigg({|v|^2\over1+|v|^2}\bigg)^{-Q/2} \bigg({|u|\over|v|}+|u|^2\bigg){1\over(1+|v|^4)^{Q/2}}dv,\,\qquad |u|<{r_m\over2}.\eqno(119)$$
The integrand in (119) is bounded above by an integrable function on $\{|v|\ge 1\}$, hence
$$J_m\le C+\int_{r_m\le|v|\le 1}\Big(|v|^{-Q-1}|u|+|v|^{-Q}|u|^2\Big)dv=C\int_{r_m}^1 \Big(r^{-2}|u|+r^{-1}|u|^2\Big)dr\le C$$
which proves (118).

To prove (117), if $|z|^2+it=|u|^2 e^{i\theta}$
and $|z'|^2+it'=|v|^2e^{i\varphi}$, then\def\zbz{z\cdot\bar {z'}}
$${1-\zbe\over1+\zbe}={|u|^2e^{i\theta}+|v|^2e^{-i\varphi}-2\zbz\over 1+|u|^2|v|^2e^{i(\theta-\varphi)}+2\zbz}$$
so, since $g$ is H\"older continuous, (117) is implied by 
$$\int_{r_m\le |v|\le 1} \big|\arg(1+|u|^2|v|^2e^{i(\theta-\varphi)}+2\zbz)\big|^{\sigma}\,|u|^{-Q}du\le C$$
and
$$\int_{r_m\le |v|\le 1} \big|\arg(|u|^2e^{i\theta}+|v|^2e^{-i\varphi}-2\zbz)+\varphi\big|^{\sigma}\,|u|^{-Q}du\le C$$
whenever  $|u|<{1\over2} r_m$. Both  these estimates follow easily as above, from the simple observation that $\arg(e^{-i\varphi}+t \omega)=-\varphi+O(t)$, as $t\to0$, if $|\omega|\le C$, uniformly in $\varphi$ (recall that $-\pi/2\le\varphi\le \pi/2$). This concludes the proof of (116) and the sharpness statement of Theorem  18.
\endpf

\vskip1.5em 


\centerline{\bf References}\bigskip
\item{[Ad1]} Adams D.R. {\sl
A sharp inequality of J. Moser for higher order derivatives},
Ann. of Math. {\bf128} (1988), no. 2, 385--398. 
\smallskip
\item{[Ad2]} Adams D.R.
{\sl Traces of potentials arising from translation invariant operators}, 
Ann. Scuola Norm. Sup. Pisa (3) {\bf 25} (1971) 203-217. 
\smallskip
\item{[Ad3]} Adams D.R., {\sl A trace inequality for generalized potentials},
Studia Math. {\bf 48} (1973), 99--105. 
\smallskip
\item{[Au1]} Aubin T., {\sl Probl\`emes isop\'erim\'etriques at espaces de Sobolev}, J. Differential Geometry {\bf11}
 (1976), 573-598.
\item{[Au2]}Aubin T.,
{\sl Meilleures constantes dans le th\'eor\`eme d'inclusion de Sobolev et un th\'eor\`eme de Fredholm non lin\'eaire pour la transformation conforme de la courbure scalaire}, J. Funct. Anal. {\bf32} (1979), 148-174. 
\smallskip
\item{[Bec]} Beckner W., {\sl Sharp Sobolev inequalities
on the sphere and the Moser-Trudin-\break ger inequality}, Ann. of
 Math. {\bf138} (1993), 213-242.\smallskip
\item{[BFM]} Branson T.P., Fontana L., Morpurgo C., {\sl Moser-Trudinger and Beckner-Onofri's inequalities on the CR sphere}, 
(2007) submitted, arXiv:0712.3905.\smallskip
\item{[BMT]} Balogh Z.M., Manfredi J.J., Tyson J.T., {\sl  Fundamental solution for the $Q$-Laplacian and sharp Moser-Trudinger inequality in Carnot groups}, J. Funct. Anal. {\bf204} (2003), 35-49.\smallskip

\item{[BCY]} Branson T.P., Chang S-Y.A., Yang P., {\sl
 Estimates and extremals for zeta function determinants on
four-manifolds}, Commun. Math.  Phys. {\bf149} (1992),
241-262.\smallskip
\item{[Ca]} Calder\'on A.-P., {\sl 
Lecture notes on pseudo-differential operators and elliptic boundary value problems, I}, 
 Consejo Nacional de Investigaciones Cientificas y Tecnicas, Instituto Argentino de Matem\'atica, Buenos Aires, 1976.\smallskip
\item {[CC]} Carleson L.,  Chang S-Y.A., {\sl On the existence of an extremal function for an inequality of J. Moser}, Bull. Sci. Math. (2) {\bf 110} (1986), 113-127.\smallskip 
\item{[Ch]} Cherrier P., {\sl Cas d'exception du theor\'eme d'inclusion de Sobolev sur les vari\'et\'es Riemanniennes et applications} Bull. Sci. Math. (2) {\bf 105} (1981) 235-288.\smallskip
\item{[CY1]} Chang S-Y.A., Yang P., {\sl Extremal             
   metrics of zeta function determinants on $4$-manifolds}, Ann. of Math.
   {\bf142} (1995),  171-212. \smallskip
\item{[CY2]} Chang S.-Y.A., Yang  P.C., {\sl  Conformal deformation of metrics on $S^2$}, J. Differential Geom. {\bf 27} (1988), 259-296.\smallskip
\item{[Ci1]} Cianchi A., {\sl Moser-Trudinger trace inequalities}, Adv. Math. {\bf 217} (2008), 2005-2044.\smallskip
\item{[Ci2]} Cianchi A., {\sl Moser-Trudinger inequalities without boundary conditions and isoperimetric problems}, Indiana Univ. Math. J. {\bf 54} (2005), 669-705. \smallskip
\item{[CoLu1]} Cohn W.S., Lu G., {\sl Best constants for Moser-Trudinger inequalities on the Heisenberg group}, Indiana Univ. Math. J. {\bf50} (2001), 1567-1591. \smallskip
\item{[CoLu2]}  Cohn W.S., Lu G., {\sl Sharp constants for Moser-Trudinger inequalities on spheres in complex space $\C^n$}, Comm. Pure Appl. Math. {\bf57} (2004), 1458-1493. \smallskip
\item{[DM]}  Djadli Z., Malchiodi A., {\sl Existence of conformal metrics with constant $Q$-curvature},  Ann. of Math. {\bf 168}  (2008), 813-858.\smallskip
\item{[F]}  Fontana L.,  {\sl Sharp borderline Sobolev inequalities on compact Riemannian manifolds}, Comment. Math. Helv. 
{\bf68} (1993), 415--454.\smallskip
\item{[FFV]} Ferone A.,  Ferone V., Volpicelli R., {\sl  Moser-type inequalities for solutions of linear elliptic equations with lower order terms}, Diff. Int. Eq. {\bf 10} (1997), 1031-1048. 
\smallskip
\item{[Fl]} Flucher M., {\sl Extremal functions for the Trudinger-Moser inequality in $2$ dimensions},  Comment. Math. Helv. {\bf 67} (1992), 471-497.
\smallskip
\item{[GT]} Gilbarg D., Trudinger N.S., {\sl Elliptic Partial Differential Equations of Second Order}, 2nd ed., Springer-Verlag, New York, 1983.\smallskip 
\item{[K]} Koldobsky A., {\sl Fourier Analysis in Convex Geometry}, Mathematical Surveys
and Monographs, vol. 116, American Mathematical Society, 2005.\smallskip 
\item{[L]} Li Y., {\sl Moser-Trudinger inequality on compact Riemannian manifolds of dimension two}, J. Part. Diff. Eq. {\bf 14} (2001), 163-192.\smallskip
\item{[LL]} Li Y., Liu P., {\sl Moser-Trudinger inequality on the boundary of compact Riemann surface}, Math. Z. {\bf 250} (2005), 363-386.\smallskip 
\item{[MMcO]} Maz'ya V., McOwen R.C., {\sl On the fundamental solution of an elliptic equation in nondivergence form}, arXiv:0806.4108v1.\smallskip
\item{[Mi]} Miranda C., {\sl Partial Differential Equations of Elliptic Type}, Springer-Verlag, New York (1970).\smallskip
\item{[Mil]} Milman E., {\sl Generalized Intersection Bodies}, J. Funct. Anal. {\bf 240} 2006, 530-567.\smallskip
\item{[Mos1]} Moser J., {\sl A sharp form of an inequality by N. Trudinger}, Indiana Univ. Math. J. {\bf20} (1970/71), 1077-1092.\smallskip
\item{[Mos2]}  Moser J., {\sl On a nonlinear problem in differential geometry}, Dynamical systems (Proc. Sympos., Univ. Bahia, Salvador, 1971), 273-280. Academic Press, New York, 1973. 
\smallskip
\item{[ON]}O'Neil R., {\sl Convolution operators in $L(p,q)$ spaces}, Duke Math. J. {\bf 30} (1963), 129-142. \smallskip
\item{[Po]} Pohozhaev  S.I., {\sl On the imbedding Sobolev theorem for pl = n}, Doklady Conference,
Section Math. Moscow Power Inst. (1965), 158-170 (Russian).\smallskip
\item{[Ruf]} Ruf B., {\sl A sharp Trudinger--Moser type inequality for unbounded domains in $\R^2$}, J. Funct. Anal. {\bf 219} (2005), 340-367.\smallskip
\item{[Tr]} Trudinger N.S., {\sl 
On imbeddings into Orlicz spaces and some applications}, 
J. Math. Mech. {\bf17} (1967), 473-483. 
\smallskip
\item{[Yu]} Yudovic V.I., {\sl Some estimates connected with integral operators and with solutions of
elliptic equations}, Dokl. Akad. Nauk SSSR {\bf138} (1961), 804-808, English translation in
Soviet Math. Doklady 2 (1961), 746-749.\smallskip
\item{[Z]} Zygmund A., {\sl Trigonometric Series}, Vol. 1, Cambridge Univ. Press (1959).

\eject
\bigskip

\noin Luigi Fontana \hskip19em Carlo Morpurgo

\noin Dipartimento di Matematica ed Applicazioni \hskip5.5em Department of Mathematics 
 
\noin Universit\'a di Milano-Bicocca\hskip 13em University of Missouri, Columbia

\noin Via Cozzi, 53 \hskip 19.3em Columbia, Missouri 65211

\noin 20125 Milano - Italy\hskip 16.6em USA 
\smallskip\noin luigi.fontana@unimib.it\hskip 15.3em morpurgoc@missouri.edu

\end